\documentclass[12pt,twoside,a4paper,leqno,onecolumn]{article}
\usepackage[english]{babel}
\usepackage[utf8]{inputenc}
\usepackage[OT1]{fontenc}
\usepackage{sectsty}
\usepackage{fncychap,fancybox}
\usepackage[print,panelleft,gray,paneltoc]{pdfscreen}
\usepackage{mathrsfs,amsthm,amsfonts,amssymb,amscd,amsbsy,amsmath}
\usepackage[retainorgcmds]{IEEEtrantools}
\usepackage{hyperref}
\numberwithin{equation}{section}
\newcommand {\R}{\mathbb{R}} 

\newtheorem{theo}{Theorem}

\newtheorem{lem}{Lemma}

\begin{document}
\pagestyle{myheadings}
\markboth{A. Hanani}{Compact hypersurfaces in a vector bundle}
\thispagestyle{empty}

\centerline{\begin{Large}\bf On compact hypersurfaces in a Riemannian\end{Large}}

\vskip2mm

\centerline{\begin{Large}\bf vector bundle with prescribed vertical\end{Large}}

\vskip2mm

\centerline{\begin{Large}\bf Gaussian curvature\end{Large}}

\vskip7mm

\centerline{\textbf{Abdellah HANANI}\footnote{Current address : Universit\'e Lille 1, UFR de Math\'ematiques, B\^at. M2, 59655, Villeneuve d'Ascq Cedex, France\\ \indent E-mail address : abdellah.hanani@math.univ-lille1.fr}}

\vskip8mm

\hrule

\vskip4mm

\noindent\textbf{Abstract. }Let M be a compact Riemannian manifold and E a Riemannian vector bundle on M. We look for hypersurfaces of E with a prescribed vertical Gaussian curvature. In trying to solve this problem fibre-wise, we loose the regularity of the resulting solution. To unsure the smoothness of the solution, we construct it as a radial graph over the unit sphere subbundle of E and prove its existence by solving in this one a nonlinear partial differential equation of Monge-Amp\`ere type.

\vskip5mm

\noindent\textbf{Keywords} : Connexions, horizontal lift, vertical lift, vertical
Gaussian curvature, Monge-Amp\`ere equations, a priori estimates, the methods.

\vskip3mm

\noindent\textbf{Mathematics Subject Classification (2010)} : 35J60, 53C55, 58G30.

\vskip4mm

\hrule

\vskip1cm

\section{Introduction}

\vskip3mm

Let $\displaystyle (M,g)$ be a compact Riemannian manifold of dimension $n\geq 1$, without boundary, and $(E,{\tilde g})$ a Riemannian vector bundle on $M$ of rank $m\geq 2$. Denote by $E_{*}$ the bundle $E$ with the zero section removed and by $\Sigma $ the corresponding unit sphere bundle. Let $VE$ be the vertical subbundle of $ TE$ and $HE$ the horizontal subbundle of $TE$ associated to a metric-connexion on $(E,{\tilde g})$. For a hypersurface ${\cal Y}$ of $E$ for which each fibre ${\cal Y}_{x}$ is a hypersurface on the fibre $E_{x}$ of $E$, the value at a point $\xi \in {\cal Y}\cap E_{x}$ of the vertical Gaussian curvature of ${\cal Y}$ is the value at the point $\xi $ of the Gaussian curvature of ${\cal Y}_{x}$ when regarded as a hypersurface of $E_{x}$.

\vskip2mm

In this study, we are interested in finding an embedding ${\cal Y}$ of $\Sigma $ into $E_{*}$ admitting a prescribed vertical Gaussian curvature equal to $K$, a given strictly positive function on $E_{*}$. We look for ${\cal Y}$ as a radial graph constructed on $\Sigma $, that is a map of the form $\xi \in \Sigma \mapsto e^{u(\xi )}\xi $, where $u\in C^{\infty }(\Sigma )$ is an unknown function extended to $E_{*} $ by letting it be radially constant. When Greek indices are used, they designate vertical directions tangent to $\Sigma $, and will range from $n+1$ to $n+m-1$. The function $u$ must satisfy on $\Sigma $ the following degenerate equation of Monge-Amp\`ere type :
\begin{equation}
det\left[\left(\delta _{\alpha }^{\beta}+D_{\alpha }uD^{\beta
}u-D_{\alpha }^{\beta }u\right)\right]=(1+v_{1})^{\frac{m+1}{2}}e^{(m-1)u}K\left(e^{u}\xi \right),
\end{equation}
where $D$ stands for the Sasaki connexion of the manifold $(E,G)$, $G$ is a
Riemannian metric on $E$ for which the vertical and the horizontal distributions are
orthogonal and $v_{1}=\sum_{n+1\leq \alpha \leq n+m-1}D_{\alpha }uD^{\alpha }u$. In
studying $(1.1)$, one needs a priori estimates on covariant derivatives of $u$ till order three.

\vskip4mm

In case the ambient space is the Euclidean one, that is when $M$ is reduced to a point, the vertical Gaussian curvature of an hypersurface of $E$ is exactly its Gaussian curvature. The question was considered at first by Oliker [9] who gave sufficient conditions on the prescribed function ensuring the existence of a solution. In particular he assumes that there exist two real numbers $r_1$ and $r_2$ such that  $0<r_1\leq 1\leq r_2$ and
\begin{equation}
K(\xi )>\Vert \xi \Vert ^{(1-m)}\ \mbox{if}\ \Vert \xi\Vert <r_1;\ K(\xi )<\Vert \xi \Vert ^{(1-m)}\ \mbox{ if} \ \Vert \xi \Vert >r_2
\end{equation} 
combined with the following monotonicity assumption :
\begin{equation}
\frac{\partial \left[\rho ^{m-1}K(\rho \xi )\right]}{\partial \rho }\leq 0, \ \mbox {for\ all\ }\xi \in \Sigma .
\end{equation}
The latter gives uniqueness up to homothety. These conditions were subsequently simplified in [3] by Delano\"e. In [1] Caffarelli, Nirenberg and Spruck were interested in finding embedded hypersurfaces of $\R^m$ whose principal curvatures satisfy a prescribed relation. To a certain extent, this question is related to the Minkowski problem (see [2], [8], [10]) from which it differs by the way of parametrising, i.e. by a radial graph instead of the inverse of the Gauss map.

\vskip4mm

The problem we deal with here requires an approach different from that used in previous works because of the degeneracy of the equation. This involves only the vertical component of the Hessian of $u$, so we have no control on the horizontal component of the second fundamental form of the hypersurface we look for, it has a strictly positive vertical Gaussian curvature but need not be convex. This makes the study of the question more interesting and shows to what extent it differs from the Euclidean case. Our
first result, which is derived in an almost elementary setting, is to
clarify this remark.

\vskip6mm
 
\begin{theo} Let $K\in \mathscr{C}^{\infty }(E_{*})$ be a strictly positive function which is constant on each fibre of E. Then there exists a radial graph ${\cal Y}$ on
$\Sigma $ whose vertical Gaussian curvature is given by $K$. Moreover, if $K$ is constant ${\cal Y}$ may be convex and in case $K$ is non constant, every such graph is non convex.
\end{theo}

\vskip4mm

Our proof relies on direct computations. The hypothesis on the prescribed function $K$
means that it is the vertical lift to $E$ of a $\mathscr{C}^{\infty }$ positive function on $M$. Any such strictly positive function gives rise to a non convex hypersurface in $E$ with strictly positive vertical Gaussian curvature. The next result is to show that assumption $(1.3)$ does not assure uniqueness not even up to homothety.

\vskip6mm

\begin{theo}Let $K\in \mathscr{C}^{\infty }(E_{*})$ be a strictly positive function such that $K(\xi )=K(\Vert \xi \Vert )$ for all $\xi $. If there exists a real number $r>0$ such that $\displaystyle K(r\xi )=r^{(1-m)}$, for all  $\xi \in \Sigma $, then there exists a radial graph ${\cal Y}$ on $\Sigma $ whose vertical Gaussian curvature is given by $K$. Such a graph may be chosen to be convex. Conversely if there exists a radial graph ${\cal Y}$ on $\Sigma $ with vertical Gaussian curvature given by $K$, then there exists a real number $r>0$ such that $\displaystyle K(r\xi )=r^{(1-m)}$, for all $\xi \in \Sigma $. Furthermore if $\displaystyle K(r\xi )=r^{(1-m)}$, for all $r>0$ and  $\xi \in \Sigma $, there exists an infinite number of non homothetic radial graphs ${\cal Y}$ on $\Sigma $ with vertical Gaussian curvature given by $K$ but a unique, up to homothety, convex one.
\end{theo}

\vskip6mm

To deal with the general case, for $u\in \mathscr{C}^{\infty }(\Sigma )$, we set
$${\cal N}_{1}(u)=det\left[(\delta _{i}^{j}+D_{i}uD^{j}u-D_{i}^{j}u)_{1\leq i,j\leq
n}\right]$$
and
$${\cal N}_{2}(u)=det\left[(\delta _{\alpha }^{\beta}+D_{\alpha }uD^{\beta
}u-D_{\alpha }^{\beta }u)_{n+1\leq \alpha ,\beta \leq n+m-1}\right].$$
Applying a continuity method in the framework of $\mathscr{C}^{\infty }$ functions [4] via the Nash and Moser inverse function theorem, for the latter see R. Hamilton [6], we first prove the following existence result.

\vskip6mm

\begin{theo}Let $f\in \mathscr{C}^{\infty}(\Sigma )$ be a strictly positive function and $\lambda $ be a strictly positive real number. Then there exists a unique solution
$u\in \mathscr{C}^{\infty }(\Sigma )$ of the equations  
$$\left\{ \begin{array}{ccl}\displaystyle {\cal N}_{1}(u)&=&1\\ \\ \displaystyle 
{\cal N}_{2}(u)&=&e^{-\lambda u}f(\xi )(1+\vert D^{v}u\vert
^{2})^{\frac{m+1}{2}}.\end{array}\right.$$
Moreover, for such a solution, the matrices $\left(G_{\alpha \beta }+D_{\alpha }uD_{\beta }u-D_{\alpha \beta }u\right)_{n+1\leq \alpha ,\beta \leq n+m-1}$ and $\left(G_{ij}+D_{i}uD_{j}u-D_{ij}u\right)_{1\leq i,j\leq n}$ are positive definite.
\end{theo}

\vskip6mm

Equation $(1.1)$ does not give any information about the horizontal behaviour of the solution if there is any. To make up for this insufficiency, we use theorem 3 to assign particular values to the horizontal derivatives. On the other hand, equation $(1.1)$ is not even locally invertible; to overcome this difficulty, we apply the fixed point theorem of Nagumo [7] to prove the following result.

\vskip6mm

\begin{theo}Let $K\in \mathscr{C}^{\infty }(E_{*})$ be an everywhere strictly positive
function. Assume that there exist two real numbers $r_1$ and $r_2$ satisfying
$0<r_1\leq 1\leq r_2$ and such that inequalities $(1.2)$ hold. Then there exists
a radial graph ${\cal Y}$ on $\Sigma $ of class $\mathscr{C}^{\infty }$ whose vertical
Gaussian curvature is given by $K$ and such that $r_1\leq \Vert \xi \Vert \leq r_2$ for all $\xi \in {\cal Y}$.
\end{theo}

\vskip6mm

The rest of this article is divided into four parts. First, we recall some preliminary results, which are needed to set the equation for prescribing the vertical Gaussian curvature on the unit sphere bundle $\Sigma $. We derive this equation in the third section by pulling the expression of the vertical Gaussian curvature back from the hypersurface. In the forth part we give the a priori estimates required in proving theorems 3 and 4 and in the last one we put these a priori estimates together and prove the results.

\vskip6mm

\section{Preliminaries and notations}

\vskip4mm

\noindent \textbf{1- }Let $(M,g)$ be a Riemannian manifold of dimension $n\geq 1$ and denote by $\nabla $ its Levi-Civita connexion. Let $(E,{\tilde g})$ be a Riemannian vector bundle on $M$ of rank $m\geq 2$, $\pi $ the projection of $E$ on $M$ and $E_{*}$ the bundle $E$ with the zero section removed. Denote by ${\tilde \nabla }$ a metric-connexion on $(E,{\tilde g})$. Let $U$ be an open set of $M$ with coordinates $(x^{i})_{1\leq i\leq n}$ and over which $E$ is trivial. The open set $\pi ^{-1}(U)$ may be equipped with a coordinate system $(x^{i},y^{\alpha })$ with $1\leq i\leq n$ and
$n+1\leq \alpha \leq n+m$, where $(y^{\alpha })_{n+1\leq \alpha \leq n+m}$ are the fibre-coordinates with respect to a fixed frame $(s_{\alpha })$ of $E$ over $U$.

\vskip4mm

We denote by $\Gamma ^{k}_{ij}$, $i,j,k\in \{1,...,n\}$, the Christoffel symbols
of the connexion $\nabla $ : $\nabla _{\varepsilon _{i}}\varepsilon _{j}=\Gamma ^{k}_{ij}\varepsilon _{k}$, where $\varepsilon _{i}={\partial }/{\partial x^i}$, and by $\Gamma ^{\beta }_{i\alpha }$, $i\in \{1,...,n\}$ and $\alpha ,\beta \in \{n+1,...,n+m\}$, the Christoffel symbols of the connexion ${\tilde \nabla }$ : ${\tilde \nabla}_{\varepsilon _{i}}s_{\alpha}=\Gamma ^{\beta}_{i\alpha}s_{\beta}$. 

\vskip4mm

On the open set $\pi ^{-1}(U)$, we can then consider the following moving frame
$${\cal S}=\{e_{i},e_{\alpha }\mid i=1,...,n\ \mbox {and}\ \alpha =n+1,...,n+m\},$$
where
\begin{equation}
e_{i}=\frac{\partial }{\partial x^i}-y ^{\alpha }\Gamma ^{\beta }_{i\alpha }\frac{\partial }{\partial y^{\beta }}
\end{equation} 
is the horizontal lift of $\varepsilon _{i}$ and $\displaystyle e_{\alpha}={\partial}/{\partial y^{\alpha}}$. Now, we equip the manifold $E$ with a Riemannian structure given by the metric $G$ defined by the following :
$$G(e_{i},e_{j})=g(\varepsilon _{i},\varepsilon _{j}),\ \ G(e_{\alpha },e_{\beta })=
{\tilde g}(s_{\alpha },s_{\beta }),\ \ G(e_{i},e_{\alpha })=0$$
and introduce the connexion $D$ of Sasaki [11] which is $G$-metric and defined as follow
\begin{equation}
D_{e_{i}}e_{j}=\Gamma ^{k}_{ij}e_{k},\ D_{e_{i}}e_{\alpha
}=\Gamma ^{\beta }_{i\alpha }e_{\beta },\ D_{e_{\alpha }}e_{i}=D_{e_{\alpha
}}e_{\beta }=0.
\end{equation}
The connexion $D$ is not torsion free. In fact, if $S^{k}_{hij}$ denotes the curvature
components of ${\tilde \nabla }$, when expressed in ${\cal S}$, the only non-zero
components of the torsion ${\cal T}$ of $D$ are those of the form
$$T^{\alpha }_{ij}=-y^{\beta }S^{\alpha }_{\beta ij}.$$
Let ${\tilde g}_{\alpha \beta }={\tilde g}(s_{\alpha },s_{\beta })$. Since ${\tilde \nabla }$ is ${\tilde g}$-metric, we obtain
$${\tilde \nabla }_{\epsilon _{i}}{\tilde g}_{\alpha \beta }={\tilde g}({\tilde \nabla }_{\varepsilon _{i}}s_{\alpha },s_{\beta })+{\tilde g}(s_{\alpha },{\tilde \nabla }_{\varepsilon _{i}}s_{\beta })=\Gamma ^{\lambda }_{i\alpha }{\tilde g}_{\lambda \beta }+\Gamma ^{\lambda }_{i\beta }{\tilde g}_{\alpha \lambda}.$$

\vskip2mm

The expression in the frame ${\cal S}$ of the components of the curvature tensor 
${\cal R}$ of $D$ are given by
$$R_{dcab}=G\Big((D_{e_{a}e_{b}}-D_{e_{b}e_{a}}-D_{[e_{a},e_{b}]})e_{c},e_{d}\Big),\
R^{d}_{cab}=G^{de}R_{ecab}$$
and standard computations yield that, for all $1\leq i,j\leq n$ and $n+1\leq \alpha ,\beta ,\lambda ,\mu \leq n+m$,
\begin{equation}
R^{i}_{\alpha \ \beta \ j}=R^{\lambda }_{\alpha \ \beta \ j}=R^{i}_{\alpha \ \beta \
\mu }=R^{\lambda }_{\alpha \ \beta \ \mu }\equiv 0.
\end{equation}

\vskip2mm

Denoting by $r$ the function $r(\xi )=\Vert \xi \Vert $ and by $\nu $ the unit
radial field which is given, on the open set $\pi ^{-1}(U)$, by 
$$\nu =r^{-1}y^{\alpha }\frac{\partial }{\partial y^{\alpha }}:=r^{-1}y^{\alpha }e_{\alpha
}.$$
Let us compute $\displaystyle D_{A}\nu $ where $A=A^{i}e_{i}+A^{\alpha }e_{\alpha }$ is a vector field on $\pi ^{-1}(U)$. At first, we have
$$\begin{array}{ccl}\displaystyle
D_{A}r&=&\displaystyle D_{A}(r^2)^{1/2}=(1/2)(r^2)^{-1/2}D_{A}r^2=(1/2)r^{-1}D_{A}({\tilde g}_{\alpha \beta }y^{\alpha }y^{\beta })\\ \\ &=&\displaystyle (1/2)r^{-1}\sum _{1\leq i\leq n}A^ie_{i}.({\tilde g}_{\alpha \beta }y^{\alpha }y^{\beta })+(1/2)r^{-1}\sum _{n+1\leq \lambda \leq n+m}A^{\lambda }\frac{\partial
}{\partial y^{\lambda }}({\tilde g}_{\alpha \beta }y^{\alpha }y^{\beta })\\ \\ &=&\displaystyle (1/2)r^{-1}\sum _{1\leq i\leq n}A^ie_{i}.({\tilde g}_{\alpha \beta }y^{\alpha }y^{\beta })+r^{-1}\sum _{n+1\leq \lambda \leq n+m}{\tilde g}_{\lambda \beta }A^{\lambda }y^{\beta }. \end{array}$$
Relations $(2.1)$ and the fact that ${\tilde \nabla }$ is ${\tilde g}$-metric imply
that
$$\begin{array}{ccl}\displaystyle e_{i}.({\tilde g}_{\alpha \beta }y^{\alpha }y^{\beta
})&=&\displaystyle \left(\frac{\partial }{\partial x^i}{\tilde g}_{\alpha \beta }\right)y^{\alpha }y^{\beta }-y^{\lambda }\Gamma ^{\mu }_{i\lambda }{\tilde g}_{\alpha \beta }\frac{\partial }{\partial y^{\mu }}(y^{\alpha }y^{\beta })\\ \\ &=&\displaystyle =2\Gamma ^{\rho }_{i\alpha }{\tilde g}_{\rho \beta }y^{\alpha }y^{\beta }-2y^{\lambda
}\Gamma ^{\mu }_{i\lambda }{\tilde g}_{\alpha \mu }y^{\alpha}=0.\end{array}$$
Therefore
$$D_{A}r=r^{-1}\sum _{n+1\leq \lambda \leq n+m}{\tilde g}_{\lambda \beta }A^{\lambda
}y^{\beta }=G(A,\nu ).$$
At present, 
$$\begin{array}{ccl}\displaystyle
D_{A}\nu &=&\displaystyle -r^{-2}(D_{A}r)y^{\alpha }e_{\alpha
}+r^{-1}(D_{A}y^{\alpha
})e_{\alpha }+r^{-1}y^{\alpha }D_{A}e_{\alpha
}\\ \\ &=&\displaystyle -r^{-1}G(A,\nu )\nu+r^{-1}\sum
_{1\leq i\leq n}A^i(D_{e_{i}}y^{\alpha }e_{\alpha }+y^{\alpha }D_{e_{i}}e_{\alpha
})\\ \\ &&\displaystyle +r^{-1}\sum _{n+1\leq \lambda \leq n+m}A^{\lambda }(D_{e_{\lambda }}y^{\alpha }e_{\alpha }+y^{\alpha }D_{e_{\lambda }}e_{\alpha
}).\end{array}$$
Taking account of $(2.2)$, we get
$$\begin{array}{ccl}\displaystyle D_{A}\nu &=&\displaystyle r^{-1}\sum _{1\leq
i\leq n}A^i(-y^{\lambda }\Gamma ^{\mu }_{i\lambda }\delta ^{\alpha
}_{\lambda }e_{\alpha }+y^{\alpha }\Gamma ^{\mu }_{i\alpha }e_{\mu
}) \\ \\ && \displaystyle -r^{-1}G(A,\nu )\nu  +r^{-1}\sum _{n+1\leq \lambda \leq n+m}A^{\lambda }\delta _{\lambda }^{\alpha }e_{\alpha }\\ \\ &=&\displaystyle r^{-1}\sum _{1\leq i\leq n}A^i(-y^{\lambda }\Gamma ^{\mu }_{i\lambda }e_{\mu }+y^{\alpha }\Gamma ^{\mu }_{i\alpha }e_{\mu })\\ \\ && \displaystyle -r^{-1}G(A,\nu )\nu 
+r^{-1}\sum _{n+1\leq \lambda \leq n+m}A^{\lambda }e_{\lambda }.\end{array}$$
Therefore
\begin{equation}
D_{A}\nu =-r^{-1}G(A,\nu )\nu  +r^{-1}\sum _{n+1\leq \lambda \leq n+m}A^{\lambda}e_{\lambda }.
\end{equation}

\vskip4mm

\noindent \textbf{2- }Let $\Sigma =\{\xi \in E\mid \Vert \xi \Vert =1\}$. The restriction of $\nu $ to $\Sigma $ is normal to  $\Sigma $. So the tangent space to $\Sigma $ at $\xi \in \Sigma $ is a direct sum of the horizontal subspace $H_{\xi }E$ of $T_{\xi }E$ and the tangent space of the fibre passing through $\xi$. This allows us to fix a local orthonormal frame field tangent to $E$ of the form  
$${\cal {R}}=\{e_{i},e_{\alpha },\nu \mid i=1,...,n\mbox{ and }\alpha =n+1,...,n+m-1\},$$
where $e_{i}$ for $i\in \{1,...,n\}$ is an horizontal vector field. Without losing
generality, we can restrict ourself to the case when $e_{i}$ is the horizontal lift of the natural vector field $\displaystyle {\partial}/{\partial x^i}$ of $M$. For $\alpha =n+1,...,n+m-1$, $e_{\alpha }$ is a vertical vector field. Let
$${\cal {R}}^{*}=\{\omega ^{A}\mid A\leq n+m\}$$
be the dual coframe. In the following, we will make use of the summation convention, when letters are used as indices they range from $1$ to $n+m$ for an upper case Latin, from $1$ to $n+m-1$ for a lower case one and from $n+1$ to $n+m-1$ for a lower case
Greek. Applying $D$ to $e_{A}$, we get a $1$-form on $E$ with values in $TE$. Expressing the result in ${\cal {R}}$ leads us to introduce the matrix $(\omega ^{A}_{B})$ of $1$-forms uniquely defined by the equalities
$$De_{A}=\omega ^{B}_{A}\otimes e_{B}.$$
From $(2.4)$, it follows that
\begin{equation}
D_{\nu }\nu =0\mbox{ and }D_{e_{a}}\nu =(1-\mu _{a})r^{-1}e_{a},\mbox{ on }\Sigma _{r},
\end{equation}
where $\Sigma _{r}=\{\xi \in E\mid \Vert \xi \Vert =r\}$ and $\mu _{a}$ is a parameter
that equals $1$ if $a$ is a horizontal direction and zero if it is a vertical one. Inserting $(2.5)$ into the previous relation, we see that, on $\Sigma _{r}$,
\begin{equation}
\omega ^{a}_{n+m}=(1-\mu _{a})r^{-1}\omega ^{a}\mbox{ and }\omega ^{n+m}_{n+m}=0.
\end{equation}
On the other hand, since $D$ is a metric connexion, it follows that
$$G(D_{e_{b}}e_{a},\nu )=-G(e_{a},D_{e_{b}}\nu ).$$
Therefore, by virtue of $(2.5)$,
\begin{equation}
\omega ^{n+m}_{a}(e_{b})=-(1-\mu _{b})r^{-1}G_{ab}\mbox{ for all }a,b\leq n+m-1.
\end{equation}

\vskip2mm

For later use, let us compute the components of the curvature tensor ${\tilde R}$ of
$\Sigma $. Using Gauss equation and relation $(2.7)$ above, which gives the
components of the second fundamental form of $\Sigma $, we show that
$${\tilde R}_{dcab}=R_{dcab}+(1-\mu _{a})(1-\mu _{b})(G_{ad}G_{bc}-G_{ac}G_{bd}).$$
Therefore, we get the values of the curvature components that will be used in next
computations :
\begin{equation}
{\tilde R}^{j}_{\alpha \beta \gamma }={\tilde R}^{j}_{\alpha \beta i}={\tilde
R}^{\gamma }_{\alpha \beta i}=0,\ n+1\leq \alpha ,\beta ,\gamma \leq n+m-1\mbox{ and }\ 1\leq i,j\leq n,
\end{equation}
and
\begin{equation}
{\tilde R}^{\lambda }_{\alpha \beta \mu }=\delta ^{\lambda }_{\beta }G_{\alpha
\mu }-\delta ^{\lambda }_{\mu }G_{\alpha \beta },\ n+1\leq \alpha ,\beta
,\lambda ,\mu \leq n+m-1.
\end{equation}

\vskip3mm

\noindent \textbf{3- }Let $u\in \mathscr{C}^{2}(\Sigma )$ be a function extended to $E_{*}$ in a radially constant way. The differential of $u$ is given by
$$du=\sum _{a=1}^{n+m-1}D_{a}u\omega ^{a}.$$
The component $D_{a}u$ is homogeneous of degree $(\mu _{a}-1)$. We
also have
$$D_{ab}u=D^2u(e_{a},e_{b})=(D_{e_{a}}Du)(e_{b}).$$
Hence
$$D_{ab}u=D_{e_{a}}\Big(Du(e_{b})\Big)-Du(D_{e_{a}}e_{b})$$
and we easily check that $\displaystyle D_{ab}u$ is homogeneous of degree $(\mu _{a}+\mu _{b}-2)$. Since $u$ is radially constant, we can write 
$$D_{a\nu }u=D_{e_{a}}\Big(Du(\nu )\Big)-Du(D_{e_{a}}\nu )=-Du(D_{e_{a}}\nu ).$$
Which in view of $(2.5)$ implies that, for all $a\leq n+m-1$,
\begin{equation}
D_{a\nu }u=-(1-\mu _{a})r^{-1}D_{a}u\mbox{ on }\Sigma _{r}.
\end{equation}
Similar computations give
\begin{equation}
D_{\nu \nu }u=0,\mbox{ on }\Sigma _{r}.
\end{equation}

\vskip6mm

\section{Derivation of the equation}

\vskip4mm

In this section, we use the method of moving frames to derive the expression of the vertical Gaussian curvature for the hypersurface under consideration. Using a homogeneity argument and pulling this expression back from the hypersurface will give us the desired equation on the unit sphere bundle $\Sigma $. For this purpose, we keep all notations of the previous part. In particular, from the choice of the local orthonormal frame ${\cal R}$ it follows that 
$${\cal {R}}_{1}=\{e_{\alpha },\nu \mid \alpha =n+1,...,n+m-1\}$$
is a local orthonormal frame field tangent to fibres of $E$. Let ${\overline D}$ denotes the induced connexion on the fibre of $E$. We look for a hypersurface ${\cal Y}$ which is a radial graph over the unit sphere bundle that is an application of the form
$${\cal Y}(\xi )=e^{u(\xi )}\xi ,\mbox{ for }\xi \in \Sigma ,$$
where $u\in C^{2}(\Sigma )$ is a function extended to $E_{*}$ by letting it be radially constant.

\vskip4mm
 
The definition of $D$ implies that $D_{e_{\alpha }}\nu $ is a vertical vector field. So, using $(2.5)$, at a point of the fibre ${\cal Y}_{x}={\cal Y}\cap E_{x}$, we obtain
\begin{equation}
{\overline D}_{e_{\alpha }}\nu =D_{e_{\alpha }}\nu =e^{-u}e_{\alpha}.
\end{equation}
Therefore, we get
\begin{equation}
{\overline D}_{\alpha \nu }u=-e^{-u}{\overline D}_{\alpha }u\mbox{ for }n+1\leq \alpha \leq n+m-1.
\end{equation}
Hence, taking $(3.1)$ into account, the definition of the covariant derivative allows us to write
\begin{equation}
D_{\nu }({\overline D}_{\alpha }u)={\overline D}_{\nu \alpha }u+{\overline
D}u({\overline D}_{\nu }e_{\alpha })=0.
\end{equation}
By a reasoning analogous to the one used to prove $(3.1)$, we show that
\begin{equation}
{\overline D}_{e_{\alpha }}e_{\beta }=D_{e_{\alpha }}e_{\beta}=\sum _{n+1}^{n+m-1}\omega ^{\gamma }_{\beta }(e_{\alpha })e_{\gamma }-e^{-u}G_{\alpha \beta}\nu .
\end{equation} 

\vskip2mm

Now, at a point $x\in M$, the tangent space of the fibre ${\cal Y}_{x}$ is spanned by the vectors $\{E_{\alpha }:={\overline D}{\cal Y}(e_{\alpha })=e_{\alpha
}+e^{u}{\overline D}_{\alpha }u\nu \}$. Therefore, on ${\cal Y}_{x}$, the induced
metric is given by 
$$h_{\alpha \beta }=G(E_{\alpha },E_{\beta })=G_{\alpha \beta }+e^{2u}{\overline D}_{\alpha }u{\overline D}_{\beta }u$$
and the unit vector field
$${\tilde \nu }=f(\nu -e^{u}{\overline D}^{\alpha }ue_{\alpha }),\ f=(1+e^{2u}{\overline D}_{\alpha }u{\overline D}^{\alpha }u)^{-\frac{1}{2}}$$
is normal to ${\cal Y}_{x}$. In view of the equalities ${\overline D}_{\nu
}u=0$, ${\overline D}_{\nu }\nu =0$, from $(3.3)$ and $(3.4)$ we obtain
\begin{equation}
\begin{array}{c}\displaystyle {\overline D}_{E_{\alpha }}E_{\beta }=\sum
_{n+1}^{n+m-1}\omega ^{\gamma }_{\beta }(e_{\alpha })e_{\gamma }-e^{-u}G_{\alpha \beta
}\nu +e^{u}(e_{\alpha }.{\overline D}_{\beta }u)\nu \\ \\ \displaystyle +{\overline D}_{\beta }ue_{\alpha }+{\overline D}_{\alpha }ue_{\beta }+e^{u}{\overline D}_{\alpha }u{\overline D}_{\beta }u\nu .\end{array}
\end{equation}
But the definition of the covariant derivative gives
$$e_{\alpha }.{\overline D}_{\beta }u={\overline D}_{\alpha \beta }u+\sum
_{n+1}^{n+m-1}\omega ^{\gamma }_{\beta }(e_{\alpha }){\overline D}_{\gamma }u.$$
Reporting into $(3.5)$, we obtain
$$\begin{array}{c} {\overline D}_{E_{\alpha }}E_{\beta }=\omega ^{\gamma
}_{\beta }(e_{\alpha })E_{\gamma }+{\overline D}_{\beta }uE_{\alpha }+{\overline D}_{\alpha }uE_{\beta }\\ \\ \displaystyle +e^{-u}[-h_{\alpha \beta }+e^{2u}{\overline D}_{\alpha \beta}u]\nu .\end{array} $$
Hence 
$$G({\overline D}_{E_{\alpha }}E_{\beta },{\tilde \nu })=fe^{-u}[-h_{\alpha \beta }+e^{2u}{\overline D}_{\alpha \beta }u].$$
This gives the components of the second fundamental form of ${\cal Y}_{x}$ when
${\cal Y}_{x}$ is regarded as a hypersurface of $E_{x}$. Namely 
$$I_{\alpha \beta }=fe^{-u}(h_{\alpha \beta }-e^{2u}{\overline D}_{\alpha
\beta }u).$$
Therefore, the Gaussian curvature of ${\cal Y}_{x}$ at the point $e^{u}\xi \in
{\cal Y}_{x}$ is
\begin{equation}
{\cal G}_{x}(e^{u}\xi )=e^{-(m-1)u}\left(1+e^{2u}\vert {\overline D}u\vert ^{2}\right) ^{-\frac{m+1}{2}}det\left(\delta _{\alpha }^{\beta }+e^{2u}{\overline D}_{\alpha }u{\overline D}^{\beta }u-e^{2u}{\overline D}_{\alpha }^{\beta }u\right).
\end{equation}
The definition of the covariant derivative gives, for $\alpha ,\beta \in \{ n+1,...,n+m-1\}$,
$$D_{\alpha \beta }u=e_{\alpha }<Du,e_{\beta }>-<Du,D_{e_{\alpha }}e_{\beta }>.$$
On the other hand, since the radial derivative of $u$ is identically equal to zero, it is clear that
$$<Du,e_{\beta }>=<{\overline D}u,e_{\beta }>$$
and taking into account the definition of the connexion $D$ which implies that 
$D_{e_{\alpha }}e_{\beta }$  is a vertical vector field, we can write
$$<Du,D_{e_{\alpha }}e_{\beta }>=<{\overline D}u,{\overline D}_{e_{\alpha
}}e_{\beta}>.$$
Hence, using Gauss equation, we get $\displaystyle D_{\alpha \beta }u={\overline
D}_{\alpha \beta }u$. Inserting into $(3.6)$, we finally arrive at the following expression of the Gaussian curvature of ${\cal Y}_{x}$ at the point $e^{u}\xi
\in {\cal Y}_{x}$ :
\begin{equation}
{\cal G}_{x}(e^{u}\xi )=e^{-(m-1)u}\left(1+e^{2u}\vert Du\vert ^{2}\right) ^{-\frac{m+1}{2}}det\left(\delta _{\alpha }^{\beta }+e^{2u}D_{\alpha }uD^{\beta }u-e^{2u}D_{\alpha }^{\beta }u\right).
\end{equation}
But the value of the vertical Gaussian curvature ${\cal G}^{v}$ of the graph ${\cal Y}$
at the point $e^{u}\xi $ is defined to be
$${\cal G}^{v}(e^{u}\xi )={\cal G}_{x}(e^{u}\xi )\ {\rm if\it \ }e^u\xi \in
{\cal Y}_{x}.$$
On the other hand, taking into account the homogeneity of the covariant derivatives of $u$, we can equate their values on ${\cal Y}_{x}$ and $\Sigma _{x}$. Therefore, by
pulling back $(3.7)$, we obtain the desired expression of the vertical Gaussian curvature of the hypersurface ${\cal Y}$ by mean of its values on $\Sigma $ :
\begin{equation}
{\cal G}^{v}(e^{u}\xi )=\left(1+\vert D^{v}u\vert ^{2}\right)^{-\frac{m+1}{2}}e^{-(m-1)u}det\left( \delta _{\alpha }^{\beta}+D_{\alpha }uD^{\beta }u-D_{\alpha }^{\beta }u\right),
\end{equation}
where $\displaystyle \vert D^{v}u\vert ^{2}=D_{\alpha }uD^{\alpha }u$.

\vskip3mm

In the sequel, we denote
$${\cal N}_{1}(u)=det\left[(\delta _{i}^{j}+D_{i}uD^{j}u-D_{i }^{j}u)_{1\leq i,j\leq n}\right]$$
and
$${\cal N}_{2}(u)=det\left[(\delta _{\alpha }^{\beta}+D_{\alpha }uD^{\beta
}u-D_{\alpha }^{\beta }u)_{n+1\leq \alpha ,\beta \leq n+m-1}\right].$$
We also denote by $G'_{u}$ the covariant $2$-tensor whose components, $G'_{ab}=G'(e_{a},e_{b})$, are given by 
$$G'_{ij}=G_{ij}+D_{i}uD_{j}u-D_{ij}u\mbox { for }1\leq i,j\leq n,$$
$$\displaystyle G'_{\alpha \beta }=G_{\alpha \beta }+D_{\alpha
}uD_{\beta }u-D_{\alpha \beta }u\mbox{ for }n+1\leq \alpha ,\beta \leq
n+m-1$$
and
$$G'_{i\alpha }=G'_{\alpha i}=0\mbox{ for }1\leq i\leq n\mbox{ and }\ n+1\leq \alpha \leq n+m-1.$$
The function $u$ is said to be admissible if the tensor $G'_{u}$ is positive definite. This allows us to view $G'_{u}$ as a new Riemannian metric on $\Sigma $.

\vskip4mm

Finally, by an adapted frame to $u$ we mean a $G$-orthonormal one that diagonalises the
Riemannian tensor $G'_{u}$.

\vskip6mm

\section{A priori estimates}

\vskip4mm

\begin{lem}Any admissible function $u\in \mathscr{C}^{2}(\Sigma )$ satisfies the
following estimate
$$\max \left(\sum _{i=1}^{n}D_{i}uD^{i}u,\sum _{\alpha =n+1}^{n+m-1}D_{\alpha
}uD^{\alpha }u\right)\leq e^{2\mathrm{osc}(u)}-1.$$
\end{lem}

\vskip5mm

\noindent\textit{Proof. }Let $u\in \mathscr{C}^{2}(\Sigma )$ be an admissible function and set $w=e^{-u}$. Easy computations show that the matrices
$$\left(wG_{ij}+D_{ij}w\right)_{1\leq i,j\leq n}$$
and
$$\left(wG_{\alpha \beta }+D_{\alpha \beta }w\right)_{n+1\leq \alpha ,\beta \leq n+m-1}$$
are positive definite. At a point $X_{0}\in \Sigma $ where the function
$\displaystyle \Omega _{1}=w^{2}+\sum _{i=1}^{n}D_{i}wD^{i}w$
attains its maximum, since $D\Omega _{1}(X_{0})=0$ in a $G$-orthonormal frame that diagonalises the symmetric matrix $(D_{ij}w)$ we get for all horizontal direction, $i\in \{1,...,n\}$,
$$D_{i}w(w+D_{ii}w)=0.$$
Thus $D_{i}w(X_{0})=0$. Since, for all $X\in \Sigma $, we have $\Omega _{1}(X)\leq \Omega _{1}(X_{0})$, we see that
$$w^{2}+\sum _{i=1}^{n}D_{i}wD^{i}w\leq w^{2}(X_{0})$$
from which, we conclude, in view of the definition of $w$, that
$$\sum _{i=1}^{n}D_{i}uD^{i}u\leq e^{2\mathrm{osc}(u)}-1.$$
Arguing analogously at a point where the function $\displaystyle \Omega _{2}=w^{2}+D_{\alpha }wD^{\alpha }w$ attains its maximum, we show that
$$\sum _{\alpha =n+1}^{n+m-1}D_{\alpha }uD^{\alpha }u\leq e^{2\mathrm{osc}(u)}-1$$
which ends the proof.

\vskip6mm

\begin{lem}Let $F\in \mathscr{C}^{3}(\Sigma \times \R)$ be a strictly positive function and $u\in \mathscr{C}^{5}(\Sigma )$ be an admissible solution of the equation
\begin{equation}
{\cal N}_{1}(u){\cal N}_{2}(u)=\left(1+\vert D^{v}u\vert ^{2}\right)^{\frac{m+1}{2}}F(\xi ,u).
\end{equation}
Assume that there exists a positive real number $C_{0}$ such that
\begin{equation}
e^{-C_{0}}\leq e^{u}\leq e^{C_{0}}.
\end{equation}
Let $\displaystyle L=\{\xi \in E\mid e^{-C_{0}}\leq \Vert \xi \Vert \leq e^{C_{0}}\}$ and denote by $\bigtriangleup u=G^{ab}D_{ab}u$ the Laplacian of $u$ with respect to the metric $G$. Then there exist positive constants $C_{1}$, $C'_{1}$ and $b$ such that
\begin{equation}
0<C_{1}\leq n+m-1+\vert Du\vert ^{2}-\bigtriangleup u\leq C'_{1}
\end{equation}
and 
$$b^{-1}G\leq G'_{u}\leq bG,$$
where $\displaystyle C_{1}=(n+m-1)(\min _{L}F)^{(n+m-1)^{-1}}$. The constants $C'_{1}$ and $b$ depend on $\displaystyle \Vert u\Vert _{\infty }$, $\displaystyle \max _{L}F$, $\displaystyle \Vert F\Vert _{\mathscr{C}^{2}(L)}$ and the geometry of $(\Sigma ,G,D)$. 
\end{lem}

\vskip5mm

\noindent\textit{Proof. }The equivalence between the two metrics becomes an obvious fact once assertion $(4.3)$ of the lemma is established. So the a priori bound on the $\mathscr{C}^{2}$-norm of $u$ comes down to establishing $(4.3)$. 

\vskip4mm

Making use of equation $(4.1)$, the admissibility of $u$ and the arithmetic and geometric means inequality, we may write, in an adapted frame to $u$,
$$\begin{array}{c}n+m-1+\vert Du\vert ^{2}-\bigtriangleup
u=G^{ab}G'_{ab}=\sum _{a=1}^{n+m-1}(1+\vert D_{a}u\vert ^{2}-D_{aa}u)\\ \displaystyle \geq (n+m-1)\left[\prod _{a=1}^{n+m-1}(1+\vert D_{a}u\vert ^{2}-D_{aa}u)\right] ^{\frac{1}{n+m-1}}\\ \displaystyle =(n+m-1)\left[F(\xi ,u)(1+\vert D^{v}u\vert ^{2})^{\frac{m+1}{2}}\right]^{\frac{1}{n+m-1}}\end{array} $$
which in view of the compactness of $L$ implies that
\begin{equation}
0<C_{1}=(n+m-1)\left(\min _{L}F\right)^{\frac{1}{n+m-1}}\leq n+m-1+\vert Du\vert ^{2}-\bigtriangleup u.
\end{equation}
On the other hand, denoting $\bigtriangleup 'u=G'^{ab}D_{ab}u$, we also have
\begin{equation}
n+m-1+\bigtriangleup 'u=G'^{ab}G_{ab}+G'^{ab}D_{a}uD_{b}u>0.
\end{equation}
\indent In view of $(4.4)$, the proof reduces to establishing an a priori bound from below on $\bigtriangleup u$. For this purpose, let $b>0$ be a fixed real number such that
$$n+m-1+\vert Du\vert ^{2}\leq b$$
and consider the function
\begin{equation}
\Gamma =(b-\bigtriangleup u)\exp\left[k\vert Du\vert ^{2}+e^{l(u+C_{0})}\right],
\end{equation}
where $k,l>0$ are real numbers to be fixed below. Let $\xi \in \Sigma $ be a point where $\Gamma $ attains its maximum and suppose that
\begin{equation}
-\bigtriangleup u(\xi )\geq 1.
\end{equation}
Hence, writing $\bigtriangleup '\log (\Gamma )\leq 0$ at $\xi$, we get
\begin{equation}
\begin{array}{c}\displaystyle \frac{-\bigtriangleup '\bigtriangleup u}{b-
\bigtriangleup u}-\frac{G'^{ab}D_{a}\bigtriangleup uD_{b} \bigtriangleup
u}{(b-\bigtriangleup u)^{2}}+l^{2}e^{l(u+C_{0})}G'^{ab}D_{a}uD_{b}u\\ \\ \displaystyle
+le^{l(u+C_{0})}\bigtriangleup 'u+2kG'^{ab}D_{abc}uD^{c}u+2kG'^{ab}G^{cd }D_{ac}uD_{bd}u\leq 0.\end{array}
\end{equation}
Next, covariantly differentiating twice the equation $(4.1)$, we obtain
\begin{equation}
G'^{ab}(D_{da}uD_{b}u+D_{a}uD_{db}u-D_{dab}u)=D_{d}\log (F)+(m+1)\frac{D_{\alpha }uD_{d}^{\alpha }u}{1+\vert D^{v}u\vert ^{2}}
\end{equation}
and 
$$\begin{array}{c}\displaystyle G'^{ab}(2D_{cda}uD_{b}u+2D_{ca}uD_{db}u-D_{cdab}u)=G'^{ae}G'^{fb}K_{cef}K_{dab}+D_{cd}\log (F)\\ \\ \displaystyle +(n+1)\frac{D_{cd\alpha }uD^{\alpha }u+D_{c\alpha }uD_{d}^{\alpha}u}{1+\vert D^{v}u\vert ^{2}}-2(m+1)\frac{D_{c\alpha }uD^{\alpha}uD_{d\beta }uD^{\beta }u}{(1+\vert D^{v}u\vert ^{2})^{2}},\end{array} $$
where $K_{cab}=D_{ca}uD_{b}u+D_{a}uD_{cb}u-D_{cab}u$. Hence contracting by $G^{cd}$ we obtain :
\begin{equation}
\begin{array}{c}\displaystyle -G'^{ab}G^{cd}D_{cdab}u=-2G'^{ab}G^{cd}D_{ca}uD_{db}u-2G'^{ab}G^{cd}D_{cda}uD_{b}u\\ \\ \displaystyle +\bigtriangleup \log (F)+(m+1)\frac{G^{cd}D_{cd\alpha}uD^{\alpha }u+D_{c\alpha }uD^{c\alpha}u}{1+\vert D^{v}u\vert ^{2}}\\ \\ \displaystyle +K-2(m+1)\frac{G^{cd}D_{c \alpha }uD^{\alpha}uD_{d\beta }uD^{\beta }u}{(1+\vert D^{v}u\vert ^{2})^{2}},\end{array}
\end{equation}
where $\displaystyle K=G^{cd}G'^{ae}G'^{fb}K_{cef}K_{dab}$. The previous
expression will imply the desired expression of $\bigtriangleup '\bigtriangleup u$ by permutations of the covariant derivatives. The last gives rise to terms involving torsion and curvature components. Recall that standard computations show that
\begin{equation}
D_{ab}u=D_{ba}u+T^{e}_{ba}D_{e}u
\end{equation}
and
\begin{equation}
D_{abc}u=D_{bac}u+R^{e}_{cba}D_{e}u+T^{e}_{ba}D_{ec}u.
\end{equation}
Consequently, we deduce that
$$G^{cd}D_{cda}u=D_{a}\bigtriangleup u+G^{cd}T^{e}_{ac}D_{ed}u+G^{cd}T^{e}_{ad}D_{ce}u+G^{cd}(R^{e}_{dac}+D_{c}T^{e}_{ad})D_{e}u.$$
Inserting into $(4.10)$ yields
\begin{equation}
\begin{array}{c}\displaystyle -G'^{ab}G^{cd}D_{cdab}u=\bigtriangleup
\log (F)-2G'^{ab}D_{ca}uD^{c}_{\ b}u-2G'^{ab}D_{a}\bigtriangleup uD_{b}u\\ \\
\displaystyle +(m+1)\frac{D_{\alpha }\bigtriangleup uD^{\alpha }u+D_{c\alpha 
}uD^{c\alpha }u}{1+\vert D^{v}u\vert ^2}-2(m+1)\frac{D_{c \alpha
}uD^{\alpha }uD^{c}_{\beta }uD^{\beta }u}{(1+\vert D^{v}u\vert ^2)^2}\\ \\ \displaystyle +(m+1)\frac{R^{ec}_{\ \ \alpha c}D_{e}uD^{\alpha }u}{1+\vert D^{v}u\vert ^2}-2G'^{ab}T^{e}_{ac}(D_{e}^{\ c}u+D^{c}_{\ e}u)D_{b}u\\ \\ \displaystyle -2G'^{ab}(R^{ec}_{\ \ ac}+D^{c}T^{e}_{ac})D_{e}uD_{b}u+K.\end{array}
\end{equation} 
Furthermore, combining $(4.11)$, $(4.12)$ and the following relation
\begin{equation}
D_{abcd}u=D_{bacd}u+R^{e}_{cba}D_{ed}u+R^{e}_{dba}D_{ce}u+T^{e}_{ba}D_{ecd}u,
\end{equation}
we check that
\begin{equation}
\begin{array}{c}\displaystyle D_{cdab}u=D_{abcd}u+T^{e}_{bc}D_{aed}u+
T^{e}_{ac}D_{ebd}u+T^{e}_{bd}D_{cae}u+T^{e}_{ad}D_{ceb}u\\ \\ \displaystyle
+(R^{e}_{dbc}+D_{c}T^{e}_{bd})D_{ae}u+(R^{e}_{dac}+D_{c}T^{e}_{ad})D_{be}u \\ \\\displaystyle +(R^{e}_{bad}+D_{a}T^{e}_{bd})D_{ce}u+(R^{e}_{bac}+D_{a}T^{e}_{bc})D_{ed}u \\ \\ \displaystyle +(D_{a}R^{e}_{dbc}+D_{c}R^{e}_{bad}+D_{ca}T^{e}_{bd}+D_{c}T^{f}_{ad}T^{e}_{bf})D_{e}u.\end{array}
\end{equation}
Therefore, taking into account $(4.11)$ and $(4.12)$, equality $(4.13)$ leads to
\begin{equation}
\begin{array}{c} \displaystyle -G'^{ab}G^{cd}D_{cdab}u=\bigtriangleup
\log (F)-2G'^{ab}D_{ca}uD^{c}_{\ b}u-2G'^{ab}D_{a}\bigtriangleup uD_{b}u\\ \\
\displaystyle +(m+1)\frac{D_{\alpha }\bigtriangleup uD^{\alpha }u+D_{c\alpha 
}uD^{c\alpha }u}{1+\vert D^{v}u\vert ^{2}}-2(m+1)\frac{D_{c \alpha
}uD^{\alpha }uD^{c}_{\beta }uD^{\beta }u}{(1+\vert D^{v}u\vert ^{2})^{2}}\\ \\ \displaystyle +(m+1)\frac{R^{ec}_{\ \ \alpha c}D_{e}uD^{\alpha }u}{1+\vert D^{v}u\vert
^{2}}+K+4G'^{ab}G^{cd}T^{e}_{ac}D_{dbe}u+E_{1}+E_{2},\end{array}
\end{equation}
where
$$\begin{array}{c}\displaystyle
E_{1}= G'^{ab}G^{cd}\left[(R^{e}_{bad}+D_{a}T^{e}_{bd})D_{ce}u+(R^{e}_{bac}+D_{a}T^{e}_{bc})D_{ed}u\right]\\ \\ \displaystyle +2G'^{ab}G^{cd}T^{e}_{ac}\left[T^{f}_{be}D_{df}u+T^{f}_{db}D_{ef}u-2D_{de}uD_{b}u\right]\end{array}$$
and
$$\begin{array}{c}\displaystyle E_{2}=2G'^{ab}G^{cd}(R^{e}_{dac}+ D_{c}T^{e}_{ad}+T^{f}_{ac}T^{e}_{df})(D_{be}u-D_{b}uD_{e}u)\\ \\ \displaystyle +G'^{ab}G^{cd}
(D_{a}R^{e}_{dbc}+D_{c}R^{e}_{bad}+D_{ca}T^{e}_{bd}+D_{c}T^{f}_{ad}T^{e}_{bf})D_{e}u\\ \\ \displaystyle +G'^{ab}G^{cd}T^{e}_{ac}(R^{f}_{edb}+R^{f}_{bde}+D_{b}T^{f}_{de}+D_{e}T^{f}_{db} +2D_{d}T^{f}_{be}D_{f}u\\ \\ \displaystyle +G'^{ab}G^{cd}T^{e}_{ac}(T^{f}_{de}T^{g}_{bf}+T^{f}_{db}T^{g}_{ef})D_{g}u.\end{array}$$
In view of $(4.7)$, the relation $G'_{ab}=G_{ab}+D_{a}uD_{b}-D_{ab}u$ and the choice of the real $b$, there exist positive constants $C_{2}$ and $C_{3}$,
independent of $u$, such that
\begin{equation}
\vert E_{1}\vert \leq C_{2}(b-\bigtriangleup u)(1+G'^{ab}G_{ab})
\end{equation} 
and
\begin{equation}
\vert E_{2}\vert \leq C_{3}(1+G'^{ab}G_{ab}).
\end{equation}
On the other hand, it is clear that
$$G^{cd}D_{c\alpha }uD^{\alpha }uD_{d\beta }uD^{\beta }u\leq \vert D^{v}u\vert
^{2}G^{ab}G^{cd}D_{ac}uD_{bd}u$$
and
\begin{equation}
G^{ab}G^{cd}D_{ac}uD_{bd}u\leq (b-\bigtriangleup u)G'^{ab}G^{cd}D_{ac}uD_{bd}u.
\end{equation}
Inserting these inequalities into $(4.16)$, by virtue of $(4.5)$, $(4.7)$ and $(4.11)$, we easily deduce the existence of a positive constant, $C_4$, such that
\begin{equation}
\begin{array}{c}\displaystyle -\bigtriangleup '\bigtriangleup u\geq
K+4G'^{ab}G^{cd}T^{e}_{ac}D_{dbe}u-2(m+2)(b-\bigtriangleup u)G'^{ab}G^{cd}D_{ac}uD_{bd}u\\ \\ \displaystyle -2G'^{ab}D_{a}\bigtriangleup uD_{b}u+(m+1)(1+\vert D^{v}u\vert
^{2})^{-1}D_{\alpha}\bigtriangleup uD^{\alpha }u\\ \\ \displaystyle +\bigtriangleup \log (F)-C_{4}(b-\bigtriangleup u)(1+G'^{ab}G_{ab}).\end{array}
\end{equation} 
\indent Now, contracting $(4.9)$ by $D^{d}u$, we get
$$G'^{ab}D_{dab}uD^{d}u=2G'^{ab}D_{da}uD^{d}uD_{b}u-D_{d}\log
(F)D^{d}u-\frac{(m+1)}{2}\frac{D_{c}\vert D^{v}u\vert ^{2}D^{c}u}{1+\vert
D^{v}u\vert ^{2}}.$$
Using $(4.11)$ and $(4.12)$, we show that
\begin{equation}
\begin{array}{c}\displaystyle G'^{ab}D_{abd}uD^{d}u=G'^{ab}D_{a}\vert Du\vert
^{2}D_{b}u-D_{d}\log (F)D^{d}u+E_{3}\\ \\ \displaystyle -\frac{m+1}{2}(1+\vert
D^{v}u\vert ^{2})^{-1}D_{\alpha }\vert Du\vert ^{2}D^{\alpha }u,\end{array}
\end{equation}
where $E_{3}$ is given by
$$E_{3}=-G'^{ab}\Big[(R^{e}_{bad}+D_{a}T^{e}_{bd}+
T^{f}_{ad}T^{e}_{bf})D_{e}u+2G_{be}T^{e}_{ad}\Big]D^{d}u+2T^{a}_{ad}D^{d}u.$$
Thus, there exists a positive constant, say $C_5$, such that
\begin{equation}
\vert E_{3}\vert \leq C_{5}(1+G'^{ab}G_{ab}).
\end{equation}
Hence, combining $(4.20)$ and $(4.21)$, we obtain, by $(4.22)$, the following relation :
\begin{equation}
\begin{array}{c} \displaystyle \frac{-\bigtriangleup '\bigtriangleup
u}{b-\bigtriangleup u}+2kG'^{ab}D_{abc }uD^{c}u\geq \frac{K}{b-\bigtriangleup u}-2(m+2)G'^{ab}D_{ac}uD_{b}^{\ c}u \\ \\ \displaystyle -\frac{m+1}{1+\vert D^{v}u\vert ^{2}}\left(\frac{-D_{\lambda}\bigtriangleup u}{b-\bigtriangleup u}+kD_{\lambda }\vert Du\vert ^2\right) D^{\lambda }u+\frac{4G'^{ab}G^{cd}T^{e}_{ac}D_{dbe}u}{b-\bigtriangleup
u}\\ \\ \displaystyle +2G'^{a b}\left(\frac{-D_{a}\bigtriangleup u}{b-\bigtriangleup u}+kD_{a}\vert Du\vert ^2\right)D_{b}u+\frac{\bigtriangleup \log (F)}{b-\bigtriangleup u}-C_{6}(1+G'^{ab}G_{ab}),\end{array}
\end{equation}
where $C_6$ is a positive constant independent of $u$. But, at the point $\xi $, where $\Gamma $ attains its maximum, the gradient of the function $\Gamma $ must vanish. Taking the logarithmic derivative of $\Gamma $, we get
\begin{equation}
\frac{-D_{a}\bigtriangleup u}{b-\bigtriangleup u}+kD_{a}\vert Du\vert ^{2}=-le^{l(u+C_{0})}D_{a}u.
\end{equation} 
So that $(4.23)$ leads to the following inequality
\begin{equation}
\begin{array}{c} \displaystyle \frac{-\bigtriangleup '\bigtriangleup
u}{b-\bigtriangleup u}+2kG'^{ab}D_{abc}uD^{c}u\geq \frac{K}{b-\bigtriangleup u}-2(m+2)G'^{ab}D_{ac}uD_{b}^{\ c}u \\ \\ \displaystyle +\frac{4G'^{ab}T^{e}_{ac}D^{c}_{\ be}u}{b-\bigtriangleup u}+\frac{\bigtriangleup \log (F)}{b-\bigtriangleup u}-2le^{l(u+C_{0})}G'^{ab}D_{a}uD_{b}u-C_{6}(1+G'^{ab}G_{ab}).\end{array}
\end{equation}
Now, we expand the following term 
$$\begin{array}{c}\displaystyle K'=G^{cd}G'^{ae}G'^{fb}\left[K_{dab}+(b-\bigtriangleup
u)^{-1}D_{a}\bigtriangleup uG'_{db}-2T^{h}_{ad}G'_{hb}\right]\\ \\ \displaystyle
\times \left[K_{cef}+(b- \bigtriangleup u)^{-1}D_{e}\bigtriangleup
uG'_{cf}-2T^{h}_{ec}G'_{hf}\right].\end{array}$$
Taking into account the choice of the real $b$, $(4.11)$, $(4.12)$ and the following two relations
\begin{equation}
G'_{ab}=G_{ab}+D_{a}uD_{b}-D_{ab}u\mbox{ and }G'^{ab}G'_{bd}=\delta ^{a}_{d},
\end{equation} 
a direct computation leads to the following
$$\begin{array}{c}\displaystyle K'\leq K-(b-\bigtriangleup
u)^{-1}G'^{ab}D_{a}\bigtriangleup uD_{b}\bigtriangleup u+4G'^{ab}G^{cd}T^{e}_{ac}D_{dbe}u\\ \\ \displaystyle +2(b-\bigtriangleup u)^{-1}\Big\{\left[(1+\vert Du\vert
^{2}+\bigtriangleup u)G'^{ab}D_{a}\bigtriangleup uD_{b}u-D^{a}\bigtriangleup
uD_{a}u\right]\\ \\ \displaystyle -G'^{ab}G^{cd}D_{a}\bigtriangleup
u\left[(R^{e}_{ cbd}+D_{b}T^{e}_{cd})D_{e}u+(T^{e}_{bc}D_{ed}u+T^{e}_{cd}D_{be}u)\right]\\ \\ \displaystyle -2G'^{ab}D_{a}\bigtriangleup
u\left[T^{c}_{bc}+(T^{e}_{bc}D^{c}u-G^{cd}T^{f}_{bc}T^{e}_{fd})D_{e}u\right]\Big\}\\ \\ \displaystyle -4G'^{ab}D_{a}u(T^{e}_{bc}D^{c}_{\ e}u+T^{e}_{bc}D^{c}uD_{e}u)-4G'^{ab}T^{e}_{ba}D_{e}u\\ \\ \displaystyle +4G^{cd}G'^{ab}G'_{ef}T^{e}_{ad}T^{f}_{bc}+4G^{ab}T^{e}_{ba}D_{e}u.\end{array}$$
Thus, taking into account $(4.2)$ and the $\mathscr{C}^{1}$ a priori estimate of Lemma 1, we conclude to the existence of a positive constant $C_{7}$ independent of $u$ such that 
$$\begin{array}{c}\displaystyle K'\leq K-(b-\bigtriangleup u)^{-1}G'^{ab}D_{a}\bigtriangleup uD_{b}\bigtriangleup u+4G'^{ab}G^{cd}T^{e}_{ac}D_{dbe}u\\ \\ \displaystyle +2(b-\bigtriangleup u)^{-1}G'^{ab}D_{a}\bigtriangleup u\Big[(1+\vert Du\vert ^2+\bigtriangleup u-G^{cd}T^{e}_{cd}D_{e}u)D_{b}u\\ \\ \displaystyle
-2T^{c}_{bc}-(R^{ec}_{\ \ bc}+G^{cd}D_{b}T^{e}_{cd}+2T^{e}_{bc}D^{c}u-2G^{cd}T^{f}_{bc}T^{e}_{fd})D_{e}u\\ \\ \displaystyle -G^{cd}(T^{e}_{cd}G_{be}+T^{e}_{bc}D_{ed}u)\Big]
+C_{7}(b-\bigtriangleup u)(1+G'^{ab}G_{ab})\\ \\ \displaystyle -2(b-\bigtriangleup
u)^{-1}D_{a}\bigtriangleup u(D^{a}u-G^{cd}T^{a}_{cd}).\end{array}$$
Using $(4.24)$, $(4.26)$ and taking $(4.7)$ into account, we easily arrive at the existence of positive constants $C_8$ and $C_9$, independent of $u$, such that
$$\begin{array}{c}\displaystyle K'\leq K-(b-\bigtriangleup u)^{-1}G'^{ab}D_{a}\bigtriangleup uD_{b}\bigtriangleup u+4G'^{ab}G^{cd}T^{e}_{ac}D_{dbe}u\\ \\ \displaystyle +C_8\left[k+le^{l(u+C_0)}\right](b-\bigtriangleup u)(1+G'^{ab}D_{a}uD_{b}u)+C_9(b-\bigtriangleup u)G'^{ab}G_{ab}.\end{array}$$
In view of the positivity of $K'$, this inequality implies that 
$$\begin{array}{c}\displaystyle K+4G'^{ab}G^{cd}T^{e}_{ac}D_{dbe}u\geq (b-\bigtriangleup u)^{-1}G'^{ab}D_{a}\bigtriangleup uD_{b}\bigtriangleup u\\ \\ \displaystyle -(b-\bigtriangleup u)\left[C_8\left(k+le^{l(u+C_{0})}\right)(1+G'^{ab}D_{a}uD_{b}u)+C_{9}G'^{ab}G_{ab}\right].\end{array}$$
Combining with $(4.25)$, we obtain 
$$\begin{array}{c}\displaystyle \frac{-\bigtriangleup '\bigtriangleup u}{b-\bigtriangleup u}-\frac{G'^{ab}D_{a}\bigtriangleup uD_{b}\bigtriangleup u}{(b-\bigtriangleup u)^{2}}+2kG'^{ab}D_{abc}uD^{c}u\geq \frac{\bigtriangleup \log (F)}{b-\bigtriangleup u}\\ \\ \displaystyle -2(m+2)G'^{ab}D_{ac}uD_{b}^{\ c}u-(C_6+C_9)(1+G'^{ab}G_{ab})\\ \\ \displaystyle -(2+C_8)\left[k+le^{l(u+C_0)}\right](1+G'^{ab}D_{a}uD_{b}u).\end{array}$$
Substituting into $(4.8)$ and taking $(4.5)$ into account, we get the following
\begin{equation}
\begin{array}{c}\displaystyle 2(k-m-2)G'^{ab}G^{cd}D_{ac}uD_{bd}u+\left[le^{l(u+C_0)}-C_6-C_9\right]G'^{ab}G_{ab}\\ \\ \displaystyle +\left[l^{2}e^{l(u+C_{0})}-(2+C_8)\left(k+le^{l(u+C_{0})}\right)\right]G'^{ab}D_{a }uD_{b}u+\frac{\bigtriangleup \log (F)}{b-\bigtriangleup u}\leq C_{10},\end{array}
\end{equation}
where $C_{10}=C_6+C_9+(n+m-1)le^{2lC_0}+(2+C_8)(k+le^{2lC_0})$. Select $k=m+2$
so that the component of $\displaystyle G'^{ab}G^{cd}D_{ac}uD_{bd}u$ vanishes. Choosing
$l$ real sufficiently large, inequality $(4.27)$ yields
\begin{equation}
G'^{ab}G_{ab}\leq C_{10}-\frac{\bigtriangleup \log (F)}{b-\bigtriangleup u}.
\end{equation}
On the other hand, by virtue of estimate $(4.2)$, the development of $\displaystyle \bigtriangleup \log (F)$ shows that there exists a positive constant $C_{11}$ such that
$$\left|\frac{\bigtriangleup \log (F)}{b-\bigtriangleup u}\right|
\leq C_{11}.$$
Therefore inequality $(4.28)$ leads to the following
\begin{equation}
G'^{ab}G_{ab}\leq C_{12}:=C_{10}+C_{11}.
\end{equation}
Finally, using equation $(4.1)$, satisfied by $u$, and the arithmetic and geometric means inequality, we check that 
$$-\bigtriangleup u\leq (1+\vert D^{v}u\vert ^{2})^{\frac{m+1}{2}}F\Big(\xi ,u(\xi )\Big)\Big(\frac{G'^{ab}G_{ab}}{n+m-2}\Big)^{n+m-2}.$$
Taking into account $(4.2)$ and the $\mathscr{C}^{1}$ a priori estimate stated in Lemma 1, inequality $(4.29)$ implies that $b-\bigtriangleup u(\xi )\leq C_{13}$,
where $C_{13}$ depends on $n$, $m$, $C_{0}$, $C_{12}$ and $\displaystyle \sup_{L}F$. The definition $(4.6)$ of $\Gamma $ allows us to conclude the proof of $(4.3)$. The equivalence between the metrics $G$ and $G'_{u}$ is clear.

\vskip4mm

\begin{lem}Keeping all the notations of the previous Lemma, let $u\in \mathscr{C}^{5}(\Sigma )$ be an admissible solution of $(4.1)$ satisfying $(4.2)$. Denote 
$$\Omega ^{2}=G'^{ab}G'^{cd}G'^{ef}D_{ace}uD_{bdf}u \mbox{ and } {\tilde
F}=\log \left[(1+\vert D^{v}u\vert ^{2})^{\frac{n+1}{2}}F(\xi ,u)\right].$$

\vskip3mm

\begin{enumerate}

\item[\rm{(i)}] There exist two positive constants $k_{1}$ and $k_{2}$ such that :
\begin{equation}
\begin{array}{c}\displaystyle \bigtriangleup '\Omega ^{2}+2G'^{ab}G'^{cd}G'^{ef}D_{ace}{\tilde F}D_{bdf}u-2G'^{ab}D_{a}\Omega ^{2}D_{b}u\\ \\ \displaystyle +H^{abcdefij}D_{ace}uD_{bdf}uD_{ij}{\tilde F}\geq -k_{1}(1+\Omega ^3)\end{array}
\end{equation}
and
\begin{equation}
\bigtriangleup '\Omega ^{2 }\geq -k_{2 }(1+\Omega ^{3})-(m+1)\frac{D_{\alpha }\Omega ^2D^{\alpha }u}{1+\vert D^{v}u\vert ^2}+2G'^{ab}D_{a}\Omega ^2D_{b }u.
\end{equation}
The components of the tensor H are given by
\begin{equation}
H^{ab...ij}=G'^{aj}G'^{ib}G'^{cd}G'^{ef}+G'^{ab}G'^{cj}G'^{id}G'^{ef}+G'^{ab}G'^{cd}G'^{ej}G'^{if}.
\end{equation}
\item[\rm{(ii)}] We have : $\Vert \Omega \Vert _{\mathscr{C}^{0}(\Sigma )}<\infty $ and, for any $\alpha \in ]0,1[$, u is uniformly bounded in $\mathscr{C}^{3,\alpha }(\Sigma )$.
\end{enumerate}
\end{lem}

\vskip5mm

\noindent\textit{Proof. }\textbf{1- }Inequality $(4.31)$ is an immediate consequence of $(4.30)$, it follows by simply expanding the terms involving ${\tilde F}$. So let us show how (ii) follows from $(4.31)$. In view of Lemma 2, there exist positive constants $C_{1}$, $C'_{1}$ and $b$ such that
\begin{equation}
0<C_{1}\leq n+m-1+\vert Du\vert ^{2}-\bigtriangleup u\leq C'_{1}\ {\rm and \it }\ b^{-1}G\leq G'(u)\leq bG.
\end{equation}
Expanding the following positive term
$$ G'^{ab}G^{cd}G^{ef}\left(D_{ace}u-\frac{1}{n+m-1}G_{ce}D_{a}\bigtriangleup
u\right)\times \left(D_{bdf}u-\frac{1}{n+m-1}G_{df}D_{b}\bigtriangleup u\right)$$
and using $(4.32)$, we show that there exists a positive constant $C_2$ such that
\begin{equation}
G'^{ab}D_{a}\bigtriangleup uD_{b}\bigtriangleup u\leq C_{2}\Omega
^{2}.
\end{equation}
On the other hand, by Cauchy's inequality and $(4.33)$, formula $(4.16)$ of the proof of Lemma 2 says that there exist two positive constants $C_3$ and $C_4$ such that
\begin{equation}
-\bigtriangleup '\bigtriangleup u\geq -C_{3}+C_{4}\Omega
^{2}+(m+1)\frac{D_{\alpha }\bigtriangleup uD^{\alpha }u}{1+\vert
D^{v}u\vert ^{2}}-2G'^{ab}D_{a}\bigtriangleup uD_{b}u.
\end{equation}

\vskip2mm

\indent Set $\Gamma =\Omega -l\bigtriangleup u$, where $l>0$ is a real number. The relation 
$$\bigtriangleup '\Omega ^{2}=2\Omega \bigtriangleup '\Omega +2G'^{ab}D_{a}\Omega D_{b}\Omega ,$$
joined to $(4.31)$ and $(4.35)$ implies the existence of positive constants $C_5$ and $C_6$ such that
\begin{equation}
\begin{array}{c}\displaystyle 2\Omega \bigtriangleup '\Gamma \geq -C_5-C_6\Omega
^{3}+lC_{4}\Omega ^{3}-2G'^{ab}D_{a}\Omega D_{b}\Omega \\ \\ \displaystyle
-\frac{2(m+1)\Omega }{1+\vert D^{v}u\vert ^2}D_{\alpha }\Gamma D^{\alpha }u+4\Omega
G'^{ab}D_{a}\Gamma D_{b}u.\end{array}
\end{equation}
At a point $\xi \in \Sigma $, where $\Gamma $ attains its maximum, we have
$$\displaystyle 0\geq 2\Omega \bigtriangleup '\Gamma \mbox{ and }D_{a}\Gamma =0.$$
So that inequality $(4.36)$ allows us to write 
$$0\geq -C_5-C_6\Omega ^3+lC_4\Omega ^3-2l^2G'^{ab}D_{a}\bigtriangleup uD_{b}\bigtriangleup u.$$
Selecting $\displaystyle l=(C_4)^{-1}(1+C_5+C_6)$ and inserting inequality $(4.34)$ into this last one, we easily obtain $\displaystyle \Omega \leq \max (1,2l^2C_2)$. From this we can easily conclude $\Omega \leq C_7$. 

\vskip3mm

Now, covariantly differentiating once equation $(4.1)$, we get 
$$\bigtriangleup '(D_{A}u)=H_{A},$$
where the right side  $\displaystyle H_{A}$ involves only covariant derivatives of $u$ of order less or equal to two. Therefore $\displaystyle \Vert H_{A}\Vert _{\mathscr{C}^{0,\alpha }(\Sigma )}\leq C_9$, and one can use Schauder's inequalities to deduce that $\displaystyle \Vert Du\Vert _{\mathscr{C}^{2,\alpha }(\Sigma )}\leq C_{10}$. 

\vskip3mm

\noindent\textbf{2- }In this paragraph we describe the steps needed to establish inequality $(4.30)$. At first, we can write
\begin{equation}
\bigtriangleup '\Omega ^{2}=\displaystyle \sum _{i=1}^{5}K_{i}.
\end{equation}
The tensor $H$ is given by $(4.32)$ and the terms $\displaystyle (K_{i})_{1\leq i\leq 5}$ are defined as follow :

\vskip4mm

$\displaystyle K_{1}=2G'^{kl}G'^{ab}G'^{cd}G'^{ef}D_{klace}uD_{bdf}u$,

\vskip3mm

$\displaystyle K_{2}=2G'^{kl}G'^{ab}G'^{cd}G'^{ef}D_{kace}uD_{lbdf}u$, 

\vskip3mm
$\displaystyle K_{3}=-2G'^{kl}H^{ab...ij}\left(D_{li}uD_{j}u+D_{i}uD_{lj}u-D_{lij}u\right)\times \left(D_{kace}uD_{bdf}u+D_{ace}uD_{kbdf}u\right)$,

\vskip3mm

$\displaystyle K_{4}=-G'^{kl}H^{ab...ij}(2D_{kli}uD_{j}u+2D_{ki}uD_{lj}u-D_{klij}u)D_{ace }uD_{bdf}u$,

\vskip3mm

$\displaystyle K_{5}=-G'^{kl}D_{k}H^{ab...ij}\left(D_{li}uD_{j}u+D_{i}uD_{lj}u-D_{lij}u\right)D_{ace}uD_{bdf}u$.

\vskip4mm

\noindent Let us write $U\simeq V$ to say that $U$ and $V$ are equivalent; i.e. if there exists a universal positive constant $c$ such that: $\vert  U(u)-V(u)\vert  \leq
c(1+\Omega ^{3})$ and use the convention of summing repeated indices from $1$ to $n+m-1$.

\vskip5mm

(i) Study of $K_{5}$. Expanding $D_{k}H^{ab...ij}$ and denoting in a $G'$-orthonormal frame
$$A=D_{kab}uD_{kcd}uD_{lca}uD_{ldb}u$$
$$B=D_{kab}uD_{kcd}uD_{lab}uD_{lcd}u.$$
Using formulas $(4.11)$ and $(4.12)$ of the proof of Lemma 2, we see that
\begin{equation}
K_{5}\simeq 6A+6B.
\end{equation}

\vskip4mm

\indent (ii) Study of $K_{4}$ . At first, we see that
$$\displaystyle K_{4}\simeq G'^{kl}H^{ab...ij}D_{klij}uD_{ace}uD_{bdf}u.$$
Thus, covariantly differentiating twice the equation satisfied by $u$, we can show
that:
$$\begin{array}{c}\displaystyle K_{4}+H^{ab...ij}D_{ij}{\tilde
F}D_{ace}uD_{bdf}u+H^{ab...ij}G'^{kp}G'^{ql }D_{jkl}uD_{iqp}uD_{ace}uD_{bdf}u\simeq \\ \\ \displaystyle \simeq G'^{kl}H^{ab...ij}D_{ace}uD_{bdf}u(D_{klij}u-D_{ijkl}u).\end{array}$$
Arguing as in the proof of Lemma 2, we see that the right hand side term is equivalent to zero. On the other hand, using $(4.11)$ and $(4.12)$ of the proof of Lemma 2, the sum of the terms in $(D^{3}u)^{4}$ in the left hand side is equivalent to $\displaystyle 3B$. Therefore
\begin{equation}
K_{4}+H^{ab...ij}D_{ij}{\tilde F}D_{ace}uD_{bdf}u+3B\simeq 0.
\end{equation}

\vskip4mm

\indent (iii) Study of $K_{3}$. By virtue of formulas $(4.11)$ and $(4.12)$, and denoting in a $G'$- orthonormal frame
$$C=D_{lab}uD_{lcd}uD_{abcd}u,$$
we see that
\begin{equation}
K_{3}\simeq 12C+2P_{abcd}D_{abcd}u,
\end{equation}
where $P_{abcd}$ are linear combinations of the derivatives of $u$ of order three, its components depend upon the curvature and torsion components as well as the covariant derivatives of $u$ of order less than two.

\vskip4mm

(iv) Study of $K_{1}$. Covariantly differentiating three times the equation satisfied by $u$, we eliminate the order five covariant derivatives of $u$. In fact, we find
$$\begin{array}{c}\displaystyle G'^{kl}G'^{ab}G'^{cd}G'^{ef}D_{acekl}uD_{bdf}u\simeq
2G'^{kl }G'^{ab}G'^{cd}G'^{ef}D_{acek}uD_{l}uD_{bdf}u\\ \\ \displaystyle
+G'^{ab}G'^{cd}G'^{ef }G'^{kj}G'^{il}\left(D_{ai}uD_{j}u+D_{i}uD_{aj}u-D_{aij}u\right)D_{cekl}uD_{bdf}u\\ \\ \displaystyle +G'^{ab}G'^{cd}G'^{ef}G'^{kj}G'^{il}\left(D_{ci}uD_{j}u+D_{i}uD_{cj}u-D_{cij}u\right)D_{aekl}uD_{bdf}u\\ \\ \displaystyle +G'^{ab }G'^{cd}G'^{ef}G'^{kj}G'^{il}\left(D_{ek}uD_{l }u+D_{k}uD_{el}u-D_{ekl}u\right)D_{acij}uD_{bdf}u\\ \\ \displaystyle -G'^{ab}G'^{cd }G'^{ef}D_{a}(G'^{kj}G'^{il})D_{ekl}uD_{cij}uD_{bdf}u\\ \\ \displaystyle -G'^{ab }G'^{cd}G'^{ef}D_{ace}{\tilde F}D_{bdf}u.\end{array}$$
Thus, a permutation of the order of covariant derivatives of $u$ allows us to write :
$$\begin{array}{c}\displaystyle G'^{kl}G'^{ab}G'^{cd}G'^{ef}D_{klace}uD_{bdf}u\simeq G'^{kl }D_{k}\Omega ^{2}D_{l}u\\ \\ \displaystyle -G'^{ab}G'^{cd}G'^{ef}D_{a}(G'^{kj}G'^{il}) D_{ekl}uD_{cij}uD_{bdf}u\\ \\ \displaystyle -3C+G'^{kl}G'^{ab}G'^{cd}G'^{ef}Q_{kace}D_{lbdf}u\\ \\ \displaystyle -G'^{ab}G'^{cd}G'^{ef}D_{ace}{\tilde F}D_{bdf}u,\end{array}$$
where $Q_{kace}$ are linear combinations of the derivatives of $u$ of order
three, its components depend upon the curvature and torsion components as well as the covariant derivatives of $u$ of order less or equal to two.

\vskip2mm

Expanding $\displaystyle D_{a}(G'^{kj}G'^{il})$ and using $(4.11)$ and $(4.12)$, we finally arrive at
\begin{equation}
\begin{array}{c}\displaystyle K_{1}\simeq 2G'^{ab}D_{a}\Omega
^{2}D_{b}u+2G'^{kl}G'^{ab}G'^{cd}G'^{ef}Q_{kace}D_{lbdf}u\\ \\ \displaystyle
-4A-6C-2G'^{ab}G'^{cd}G'^{ef}D_{ace}{\tilde F}D_{bdf}u.\end{array}
\end{equation}

\indent (v) Inserting $(4.38)$, $(4.39)$, $(4.30)$ and $(4.31)$ into $(4.37)$, we obtain
\begin{equation}
\begin{array}{c}\displaystyle \bigtriangleup '\Omega ^{2}\simeq 2G'^{ab}D_{a}\Omega ^{2}D_{b}u+2(P+Q)_{abcd}D_{abcd}u+2A\\ \\ \displaystyle +3B+6C+K_{2}-H^{ab...ij}D_{ij}{\tilde F}D_{ace}uD_{bdf}u\\ \\ \displaystyle -2G'^{ab}G'^{cd}G'^{ef}D_{ace}{\tilde F}D_{bdf}u.\end{array}
\end{equation}
Let us introduce the tensors $U$ and $V$ which are defined in a $G'$-orthonormal
frame by
$$U_{abcd}=D_{abcd}u+D_{lab}uD_{lcd}u+(P+Q)_{abcd}$$
and 
$$V_{abcd}=D_{abcd}u+D_{lab}uD_{lcd}u+D_{lac}uD_{lbd}u.$$
We easily check that
$$\Vert U\Vert ^2\simeq \frac{1}{2}K_{2}+2C+B+2(P+Q)_{abcd}D_{abcd}u$$
and, using once more relation $(4.12)$ of the proof of Lemma 2, we obtain
$$\Vert V\Vert ^2\simeq \frac{1}{2}K_{2}+2A+2B+4C.$$ 
Thus, $(4.42)$ becomes
$$\begin{array}{c}\displaystyle \bigtriangleup '\Omega ^2\simeq 2G'^{ab}D_{a}\Omega
^2D_{b}u-H^{ab...ij}D_{ij}{\tilde F}D_{ace}uD_{bdf}u\\ \\ \displaystyle +\Vert U\Vert ^2+\Vert V\Vert ^2-2G'^{ab}G'^{cd }G'^{ef}D_{ace}{\tilde F}D_{bdf}u\end{array}$$
and $(4.30)$ follows from this relation in view of the positivity of $\Vert U\Vert
^2$ and $\Vert V\Vert ^2$.

\vskip6mm

\section{Proof of the results}
\subsection{Proof of theorem 1}

\vskip4mm

In order to end the proof of theorem 1, let us consider the function 
$$u=\frac{-\log (K)}{m-1}.$$
From the assumption on $K$ we see that $u=u\circ \pi $ in $E_{*}$. Let $\pi _{*}$ stand for the tangent map of $\pi $. From the usual rules of differentiation it follows that
$$Du(Y)=Du\left( \pi _{*}(Y)\right) \mbox{ for }Y\in TE.$$
Recall that for any $\xi \in E_{*}$ the vertical subspace $V_{\xi }E$ of $T_{\xi }E$ is exactly the kernel of $\pi _{*\mid T_{\xi }E}$. Thus, for all vertical direction $\alpha $, we get $\displaystyle D_{\alpha }u=0$ and, since $D_{e_{\alpha }}e_{\beta }$ is a vertical vector field for $\alpha ,\ \beta \in \{n+1,...,n+m-1\}$, the definition of the covariant derivative implies that $D_{\alpha \beta }u=0$. Inserting this into equation $(3.8)$, we deduce that 
$${\cal G}^{v}(e^{u}\xi )=e^{-(m-1)u(\xi )}=K(\xi )$$
as claimed.

\vskip4mm

Now, if $K$ is constant, then the radial graph ${\cal Y}:\xi \mapsto e^{u(\xi )}\xi $,
with $u$ as above, is a convex hypersurface whose vertical Gaussian curvature is given by $K$. Conversely, assume that there exists a function $u\in \mathscr{C}^{\infty }(\Sigma )$ such that
\begin{equation}
det\left( \delta _{\alpha }^{\beta}+D_{\alpha }uD^{\beta }u-D_{\alpha }^{\beta
}u\right)=\left( 1+\vert D^{v}u\vert ^{2}\right) ^{\frac{m+1}{2}}e^{(m-1)u}K(\xi )
\end{equation}
and that the hypersurface ${\cal Y}$ given by 
$${\cal Y}(\xi )=e^{u(\xi )}\xi ,\mbox{ for }\xi \in \Sigma ,$$
is convex so that its second fundamental tensor is a $2$-covariant positive one. We
need to compute the components of this tensor. Keeping all previous notations and following the formalism of moving frame we check that these are given by the $(n+m-1)\times (n+m-1)$-matrix
$$g'=\Big((1-\mu _{a})G _{ab}+(1-2\mu _{a}){\tilde D}_{a}u{\tilde
D}_{b}u-{\tilde D}_{ab}u\Big) _{1\leq a,b\leq n+m-1},$$
where ${\tilde D}_{a}u=e^{\mu _{a}u} D_{a}u$ and ${\tilde D}_{ab}u=e^{(\mu
_{a}+\mu _{b})u}D_{ab}u$. By the assumption on ${\cal Y}$, the matrix $g'$ is
positive. Particularly, for any real $\varepsilon >0$, the $n\times n$-matrix
\begin{equation}
\left(\varepsilon G_{ab}-D_{a}uD_{b}u-D_{ab}u\right)_{1\leq a,b\leq n}\mbox{ is
positive definite.}
\end{equation}
Now, consider the function
$$\Gamma (u)=-\epsilon u+G^{ab}(\mu_{a}D_{a}u)(\mu_{b}D_{b}u)$$
at a point $\xi _{0}\in \Sigma $ where it attains its maximum. For all horizontal direction $i\in \{ 1,...,n\} $, we get $D_{i}\Gamma =0$ so that, in an adapted frame to $u$, we can write
$$D_{i}u(-\epsilon +D_{ii}u)=0.$$
Taking $(5.2)$ into account, for all horizontal direction $i$, we obtain $D_{i}u(\xi
_{0})=0$. Consequently $\Gamma (u)(\xi _{0})=-\varepsilon u(\xi _{0})$ and since $\Gamma (u)$ attains its maximum at $\xi _{0}$, we conclude that there exists a positive constant $C$ such that, for any $\varepsilon >0$,
$$G^{ab}(\mu_{a}D_{a}u)(\mu_{b}D_{b}u)\leq C\varepsilon .$$
From this we conclude that the horizontal component of the gradient of $u$ is identically equal to zero. Using the definition of the covariant derivative, we show that, for any horizontal direction $i$, $1\leq i\leq n$, and any vertical direction $\alpha $, $n+1\leq \alpha \leq n+m-1$,
\begin{equation}
D_{i\alpha }u=D_{\alpha i}u=0.
\end{equation}
The first equality in $(5.3)$ follows from the fact that torsion components of the form
$T^{a}_{i\alpha }$ are zero. Taking into account relation $(2.8)$, we see that $\displaystyle D_{i\alpha \beta }u=D_{\alpha \beta i}u$. Now, for any vertical directions $\alpha ,\beta $ and any horizontal one $i$, the definition of the covariant derivative yields
$$D_{i\alpha \beta }u=0.$$
Therefore, covariantly differentiating the equation $(5.1)$ in any horizontal direction, we see that the horizontal gradient of $K$ is identically equal to zero and as it is constant on each fibre of $E$, the vertical component of its gradient is also identically equal to zero. We then conclude that $K$ is constant.

\vskip6mm

\subsection{Proof of theorem 2}

\vskip4mm

Relying crucially on computations from section 2, we deduce that the vertical Gaussian curvature of the bundle $\Sigma _{r}$ is $\displaystyle r^{-(m-1)}$. The bundle $\Sigma _{r}$ is closed and convex. Now suppose that there exist a radial graph over $\Sigma $ with vertical Gaussian curvature given by $K$. Then there exists $u\in \mathscr{C}^{\infty }(\Sigma )$ satisfying $(5.1)$. Set
$$\psi (r)=r^{m-1}K(r),\ R_{1}=\min _{\Sigma }e^{u}\mbox{ and }R_{2}=\max _{\Sigma }e^{u}.$$
From $(5.1)$, we see that
$$\psi (R_{1})\leq 1\mbox{ and }\psi (R_{2})\geq 1$$
thus, by continuity, there exists a real $r>0$ such that $\psi (r)=1$.

\vskip3mm

Now suppose that, for any $r>0$ and any $\xi \in \Sigma $,
\begin{equation}
K(r\xi )=K(r)=\frac{1}{r^{m-1}}.
\end{equation} 
If $u\in \mathscr{C}^{\infty }(\Sigma )$ satisfies $(1)$, then for any non constant $v\in \mathscr{C}^{\infty }(M)$, the function ${\tilde u}=u+v\circ \pi $ is also a solution of $(5.1)$ and clearly their graphs are not homothetic. Thus, there is an infinite number of non homothetic closed hypersurfaces with vertical Gaussian curvature given by $K$. But two convex closed hypersurfaces with vertical Gaussian curvature given by $K$ are necessarily homothetic. In fact, if ${\cal Y}$, the radial graph of a function $u\in \mathscr{C}^{\infty }(\Sigma )$, is such an hypersurface, then
\begin{equation}
det\left( \delta _{\alpha }^{\beta}+D_{\alpha }uD^{\beta }u-D_{\alpha }^{\beta
}u\right)=\left( 1+\vert D^{v}u\vert ^{2}\right) ^\frac{{m+1}{2}}.
\end{equation}
Recall from the proof of theorem 2 that the convexity of the hypersurface ${\cal Y}$ implies the nullity of the horizontal gradient of $u$. On the other hand, at a point
$\xi \in \Sigma $ where $u$ attains its maximum, and in a $G$-orthonormal frame in
which the symmetric matrix $\displaystyle \left[D_{\alpha }uD_{\beta }u-D_{\alpha \beta }u\right]$ is diagonal, it is clear that the eigenvalues of the matrix
$$\left(G'_{u}=G_{\alpha \beta }+D_{\alpha }uD_{\beta }u-D_{\alpha \beta }u\right)_{\alpha
\beta }$$
are strictly positive at $\xi $. By continuity, $\displaystyle \left[(G'_{u})_{\alpha \beta}\right]$ must be positive definite everywhere. Let us introduce the operator $\bigtriangleup '$ by setting $\bigtriangleup 'u=G'^{\alpha \beta }D_{\alpha \beta
}u$. Covariantly differentiating $(5.5)$, we get 
$$G'^{\alpha \beta }(D_{\lambda \alpha }uD_{\beta }u+D_{\alpha }uD_{\lambda
\beta }u-D_{\lambda \alpha \beta }u)=\frac{m+1}{2}\frac{D_{\lambda }\vert
D^{v}u\vert ^2}{1+\vert D^{v}u\vert ^2}$$
and contracting by $D^{\lambda }u$, we obtain
\begin{equation}
G'^{\alpha \beta }D_{\lambda \alpha \beta }uD^{\lambda }u=G'^{\alpha \beta
}D_{\alpha }\vert D^{v}u\vert ^2D_{\beta }u-\frac{(m+1)}{2}\frac{D_{\lambda
}\vert D^{v}u\vert ^2}{1+\vert D^{v}u\vert ^2}.
\end{equation}
From the expressions $(2.8)$ and $(2.9)$ of the curvature components, we can write 
$$G'^{\alpha \beta }D_{\lambda \alpha \beta }uD^{\lambda }u=G'^{\alpha \beta
}D_{\alpha \beta \lambda }uD^{\lambda }u+G'^{\alpha \beta }D_{\alpha
}uD_{\beta }u-\vert D^{v}u\vert ^{2}G'^{\alpha \beta }G_{\alpha \beta
}.$$
Inserting into $(5.6)$, and evaluating at a point $\xi \in \Sigma $ where the function
$$\Gamma =\frac{1}{2}G^{\alpha \beta }D_{\alpha }uD_{\beta }u$$
attains its maximum, we get
$$G'^{\alpha \beta }D_{\alpha \beta \lambda }uD^{\lambda }u=\vert D^{v}u\vert
^{2}G'^{\alpha \beta }G_{\alpha \beta }-G'^{\alpha \beta }D_{\alpha }uD_{\beta }u.$$
But, at the point $\xi $, we must have $\bigtriangleup '\Gamma \leq 0$. Thus
\begin{equation}
\vert D^{v}u\vert ^2G'^{\alpha \beta }G_{\alpha \beta }-G'^{\alpha \beta
}D_{\alpha }uD_{\beta }u+G'^{\alpha \beta }G^{\lambda \mu }D_{\alpha \lambda }uD_{\beta \mu }u\leq 0.
\end{equation}
In a $G$-orthonormal frame in which the symmetric matrix $\displaystyle (D_{\alpha }uD_{\beta }u-D_{\alpha \beta }u)$ is diagonal, we can write
$$G'^{\alpha \beta }D_{\alpha }uD_{\beta }u=\sum _{\alpha =n+1}^{n+m-1}\frac{\vert
D_{\alpha }u\vert ^{2}}{1+\vert D_{\alpha }u\vert ^{2}-D_{\alpha \alpha }u}\leq
\vert D^{v}u\vert ^{2}G'^{\alpha \beta }G_{\alpha \beta }.$$
Therefore, the relation $(5.7)$ implies
$$0\leq G'^{\alpha \beta }G^{\lambda \mu }D_{\alpha \lambda }uD_{\beta \mu }u\leq 0$$
and in particular, at the point $\xi $, we must have for any directions $\alpha $ and
$\beta $,
$$D_{\alpha \beta }u=0.$$
Inserting into $(5.5)$, we obtain
$$1+\vert D^{v}u\vert ^{2}= det(\delta _{\alpha }^{\beta}+D_{\alpha }uD^{\beta
}u)=(1+\vert D^{v}u\vert ^{2})^{\frac{m+1}{2}}.$$
Hence $\vert D^{v}u\vert ^{2}(\xi )=0$ and then $\vert D^{v}u\vert ^{2}=0$
everywhere. Thus, taking into account the fact that the horizontal gradient of $u$ is
identically equal to zero, the function $u$ must be a constant. This ends the proof of the theorem.

\vskip7mm

\subsection{Proof of theorem 3} 

\vskip4mm

Let $f\in \mathscr{C}^{\infty}(\Sigma)$ be a strictly positive function and $\lambda >0$ a real number. We want to solve the equations
\begin{equation}
\left\{ \begin{array}{c}\displaystyle {\cal
N}_{1}(u)=1\\ \\ \displaystyle 
{\cal N}_{2}(u)=e^{-\lambda u}f(\xi
)(1+\vert D^{v}u\vert ^{2})^{\frac{m+1}{2}}.\end{array}\right.
\end{equation}
Remark first that a solution of $(5.8)$, if there is any, is necessarily admissible. In fact, at a point  $\xi \in \Sigma $ where $u$ attains its maximum, and in a frame adapted to $u$, it is clear that the eigenvalues of the tensor $G'_{u}$ are strictly positive at $\xi$. By continuity and since $\displaystyle {\cal N}_{1}(u)>0$ we see that the matrix $\displaystyle \left[(G'_{u})_{ij}\right]_{1\leq i,j\leq n}$ must be
positive definite everywhere. We also have ${\cal N}_{2}(u)>0$ so that the matrix $\displaystyle \left[(G'_{u})_{\alpha \beta }\right]_{n+1\leq \alpha ,\beta \leq n+m-1}$ is positive definite everywhere.

\vskip3mm

On the other hand, if a solution exists, it is unique for if there exist two solutions $u_{1}$ and $u_{2}$, we put $u=u_{1}-u_{2}$ and let $\xi \in \Sigma $ be a point 
where $u$ attains its minimum. In an orthonormal frame for $\displaystyle
[(G'_{u_{1}})_{\alpha \beta }]$ that diagonalises $\displaystyle [(G'_{u_{2}})_{\alpha \beta }]$ we have $(G'_{u_{2}})_{\alpha \beta }=(G'_{u_{1}})_{\alpha \beta }+D_{\alpha
\beta}u$. Thus, we may write :
$$\frac{{\cal N}_{2}(u_{2})}{{\cal N}_{2}(u_{1})}=\prod _{\alpha =n+1}^{n+m-1}(1+D_{\alpha \alpha }u)=e^{\lambda u}.$$
Taking into account the fact that $D_{\alpha \alpha }u(\xi )\geq 0$, the last equality
implies : $u(\xi )\geq 0$. Hence $u_{1}-u_{2}\geq 0$ everywhere. By an analogous
argument, but at a point where $u$ attains its maximum, we get the inequality in the other sense. Thus, the functions $u_{1}$ and $u_{2}$ must be equal.

\vskip3mm

To treat the existence part, we will use the continuity method in the framework of 
$\mathscr{C}^{\infty }$ functions [4]. This method consists of the following steps.  

\vskip2mm

For $t\in [0,1]$, consider the family of equations
\begin{equation}
\left\{ \begin{array}{c}\displaystyle {\cal N}_{1}(u)=1\\ \\ \displaystyle 
{\cal N}_{2}(u)=e^{-\lambda u}\Big[f(\xi )(1+\vert D^{v}u\vert ^{2})^{\frac{m+1}{2}}\Big]^{t}.\end{array}\right.
\end{equation}
By the previous discussion, any solution $u_{t}$ of $(5.9)$ is admissible. So let
$T$ be the set of $t\in [0,1]$ for which such a solution exists in $\mathscr{C}^{\infty}(\Sigma )$.

\vskip2mm

Observe that the function $u_{0}=0$ is a solution of $(5.9)$ for $t=0$. Hence $T$ is not empty. If we prove that $T$ is open and closed in $[0,1]$ then $T=[0,1]$ and the
equation $(5.8)$ is solved in $\mathscr{C}^{\infty }(\Sigma)$. Recall that the $\mathscr{C}^{\infty }$ regularity follows from the well known regularity theory of elliptic equations.

\vskip2mm

For the closeness, we have to establish, for any $k\geq 0$, a uniformly bound
of the $\mathscr{C}^{k}$-norm of any admissible solution of $(5.9)$. The $\mathscr{C}^{0}$ estimate $\displaystyle \vert u_{t}\vert \leq \lambda ^{-1}\Vert \log (f)\Vert _{\infty }$ of any solution of $(5.9)$ is immediate, for if $\xi _{1}\in \Sigma $ is a point where a solution $u_{t}$ of $(5.9)$ attains its maximum, then in a frame adapted to $u_t$ we get, for any $\alpha $, $D_{\alpha }u_{t}(\xi
_{1})=0$ and $D_{\alpha \alpha}u_{t}(\xi _{1})\leq 0$. So $(5.9)$ implies :
$$0\leq \log [{\cal {N}}_{2}(u_{t})](\xi _{1})=-\lambda u_{t}+t\log [f(\xi _{1})].$$
Hence $\displaystyle u_{t}\leq \lambda ^{-1}\Vert \log (f)\Vert _{\infty }$. We
complete our $\mathscr{C}^{0}$ estimate by an analogous reasoning, considering a point where $u_t$ attains its minimum.

\vskip2mm

The first, second and third order a priori estimates are given in Lemmas 1, 2 and 3 and the higher order estimates are established by induction; we apply the maximum principle as in lemma 3 to functional depending on a norm of derivatives of order $k$ we want to bound and on already bounded quantities. These estimates may be recovered by Schauder's inequalities applied inductively to equations obtained by covariantly differentiating the initial one. In fact, it follows from the $\mathscr{C}^3$-estimates that the components of the linearised operator are bounded in $\mathscr{C}^{0,\alpha }(\Sigma )$. Thus keeping the same notations as in lemma 3 and assuming that 
$$\Vert u\Vert _{\mathscr{C}^{k,\alpha }(\Sigma )}\leq C_{k},\ k\geq 3,$$
where $u\in \mathscr{C}^{\infty }(\Sigma )$ is an admissible solution of the equation
$$\displaystyle {\cal N}_{1}(u){\cal N}_{2}(u)=F(\xi
,Du,u).$$
Covariantly differentiating $(k-1)$ times this equation, on any open set $U$ with
coordinates, we get 
$$\bigtriangleup '(D_{a_{1}...a_{k-1}}u)=H_{a_{1}...a_{k-1}},$$
where the right side $\displaystyle H_{a_{1}...a_{k-1}}$ involves only covariant derivatives of $u$ of order less or equal to $k$. Therefore $\displaystyle \Vert
H_{a_{1}...a_{k-1}}\Vert _{\mathscr{C}^{0,\alpha }(\Sigma )}\leq Cste$. We can apply Schauder's inequalities to deduce that for any compact $K\subset U$, we have $\displaystyle \Vert D^{(k-1)}u\Vert _{\mathscr{C}^{2,\alpha }(K)}\leq C'_{k+1}$. Consequently, $\displaystyle \Vert u\Vert _{\mathscr{C}^{k+1,\alpha }(\Sigma )}\leq C_{k+1}$.

\vskip3mm

To show that $T$ is open, let $A^{\infty }(\Sigma)$ be the set of admissible functions
$u\in \mathscr{C}^{\infty }(\Sigma )$ and denote by $\Theta$ the following subset : 
$$\Theta =\{u\in A^{\infty }(\Sigma )\mid {\cal N}_{1}(u)=1\}.$$
The set $\Theta $ is a hypersurface of $A^{\infty }(\Sigma )$, this will be shown below. Now let $\Gamma $ be the functional defined on $\Theta \times [0,1]$ by
$$\Gamma (u,t)=\log \left[{\cal N}_{2}(u)\right]+\lambda u-t\log \left[f(\xi )(1+\vert D^{v}u\vert ^2)^{\frac{m+1}{2}}\right].$$
The function $\Gamma $ is continuously differentiable and its differential at $u\in
\Theta $ is a linear operator from $B^{\infty }(\Sigma )$ into $\mathscr{C}^{\infty
}(\Sigma )$ given by
$$\displaystyle d_{u}\Gamma (w)=G'^{\alpha \beta }(2D_{\alpha }uD_{\beta }w-D_{\alpha \beta }w)+\lambda w-t(m+1)\frac{D^{\alpha }uD_{\alpha }w}{1+\vert D^{v}u\vert ^{2}},$$
where 
$$B^{\infty }(\Sigma )=\Big\{w\in \mathscr{C}^{\infty }(\Sigma )\mid \sum _{1\leq
i,j\leq n}G'^{ij}(2D_{i}uD_{j}w-D_{ij}w)=0\Big\}.$$
For $u$ given in $\Theta $, $d_{u}\Gamma $ is invertible. In fact, its null space is trivial since $\lambda >0$. So we are done if we prove that, for any $v\in \mathscr{C}^{\infty}(\Sigma )$, there exists a solution $w\in B^{\infty }(\Sigma )$ of the equation
\begin{equation}
d_{u}\Gamma (w)=v.
\end{equation}
The proof is the same for arbitrary values of the parameter $t$. So assume $t=0$, set $r$ for the function $r(\xi )=\Vert \xi \Vert $ for all $\xi \in E_{*}$, extend $u$ and $v$ to $E_{*}$ as radially constant functions and set $\displaystyle {\tilde \Sigma }=\{\xi \in E\mid 1\leq \Vert \xi \Vert \leq 2\}$. Let us also denote by $L[u]$ the linear operator defined by
$$\displaystyle L[u](w):=r\sum _{1\leq i,j\leq n}G'^{ij}_{u}D_{ij}w+r^{2}\sum _{m+1\leq \alpha ,\beta \leq n+m-1}G'^{\alpha \beta }_{r,u}D_{\alpha \beta }w-\lambda w,$$
where 
$$\displaystyle G'^{\alpha \beta }_{r,u}=G^{\alpha \beta }+r^{2}D^{\alpha }uD^{\beta
}u-r^{2}D^{\alpha \beta }u.$$
For $s\in [0,1]$ and $\displaystyle y\in \mathscr{C}^{\infty }({\tilde \Sigma })$, let $H_{s}y=w_{s}$ be the unique solution of the problem
\begin{equation}
\left\{ \begin{array}{cc}\displaystyle D_{\nu \nu }w_{s}+L[u](w_{s})-a
w_{s}=sV(y)+D_{\nu \nu }y+rBD_{\nu }y\mbox{ in }{\tilde \Sigma }\\ \\ \displaystyle
D_{\nu }w_{s}=0\mbox{ on }\partial{\tilde \Sigma },\end{array}\right.
\end{equation}
where $\nu $ is the unit radial field, $a\geq 0$ is a real number,
$$V(y)=2r\sum _{1\leq i,j\leq n}G'^{ij}_{u}D_{i}uD_{j}{\tilde y}+2r^{2}\sum _{n+1\leq \alpha ,\beta \leq n+m-1}G'^{\alpha \beta }_{r,u}D_{\alpha
}uD_{\beta }{\tilde y}-v(1+ar)$$
and 
$$B=m-1+r^{2}\sum _{n+1\leq \alpha \leq n+m-1}(D_{\alpha}uD^{\alpha}u-D_{\alpha }^{\alpha }u).$$
The function ${\tilde y}$ stands for the extension as a radially constant function to
${\tilde \Sigma }$ of the restriction to $\Sigma $ of $y$.

\vskip2mm

For a proof, we use similar arguments to those of Gilbarg and Tr\"udinger [5], theorem 6.31. Moreover, since problem $(5.11)$ is uniformly elliptic, if ${\cal B}$ is a bounded subset of $\mathscr{C}^{\infty}(\Sigma )$, there exist a sequence of positive real numbers $C_{k}$ such that
\begin{equation}
\Vert w_{s}\Vert _{\mathscr{C}^{k}({\tilde \Sigma })}\leq C_{k}\ \mbox{for\ any\ }(s,y)\in [0,1]\times {\cal B}.
\end{equation}  
It follows that the operator $H$ defined on $[0,1]\times \mathscr{C}^{\infty }({\tilde \Sigma})$ by $\displaystyle H(s,y)=y-H_{s}y$ is compact. Our aim is to solve equation 
\begin{equation}
H(s,w)=0,\mbox{ for }s\in [0,1].
\end{equation}
Any solution of $(5.13)$ is a radial constant. In fact, such a solution satisfies :
\begin{equation}
\left\{ \begin{array}{cc}\displaystyle L[u](w)-aw=sV(w)+rBD_{\nu }w\; \mbox{ in }
{\tilde \Sigma }\\ \\ \displaystyle D_{\nu }w=0\; \mbox{ on }\partial {\tilde \Sigma
}.\end{array}\right.
\end{equation}
Recall from section 2 that
\begin{equation}
D_{e_{a}}(r\nu )=(1-\mu _{a})e_{a},
\end{equation}
where $\mu _{a}$ is equal to $1$ or $0$ accordingly to whether $e_a$ is horizontal or
vertical. Hence, using the definition of covariant derivative, we can write :
$$D_{e_{a}}(rD_{\nu }w)=D^{2}w(e_{a},r\nu )+Dw\Big(D_{e_{a}}(r\nu )\Big)=rD_{\nu}(D_{a}w)+(1-\mu _{a})D_{a}w.$$
From this, one deduces that 
\begin{equation}
rD_{\nu }(D_{a}w)=D_{e_{a}}(rD_{\nu }w)-(1-\mu _{a})D_{a}w.
\end{equation}
Using again the definition of the covariant derivative, we can at first write :
$$D_{ab}(rD_{\nu }w)=D^{2}(rD_{\nu }w)(e_{a},e_{b})=e_{a}\Big[D(rD_{\nu }w)(e_{b})\Big]-D(rD_{\nu }w)(D_{e_{a}}e_{b})$$
and then
$$\begin{array}{c}\displaystyle D_{ab}(rD_{\nu }w)=e_{a}\left[D_{e_{b}}\left(Dw(r\nu
)\right)\right]-\left(D_{e_{a}}e_{b}\right)Dw(r\nu )\\ \\ \displaystyle
=D_{e_{a}}(D^{2}w(e_{b},r\nu )+D_{e_{a}}Dw\left(D_{e_{b}}(r\nu
)\right)\\ \\ \displaystyle -D^{2}w\left(D_{e_{a}}e_{b},r\nu
\right)-Dw\left(D_{D_{e_{a}}e_{b}}(r\nu )\right)\end{array}.$$
Therefore
$$\begin{array}{c}\displaystyle D_{ab}(rD_{\nu }w)=D^{3}w(e_{a},e_{b},r\nu)+D^{2}w\left(e_{b},D_{e_{a}}(r\nu )\right) +D^{2}w\left(e_{a},D_{e_{b}}(r\nu )\right)\\ \\
\displaystyle +Dw\left(D_{e_{a}}(D_{e_{b}}(r\nu ))\right)-Dw\left(D_{D_{e_{a}}e_{b}}(r\nu )\right).\end{array}$$
But relation $(5.15)$ implies that
$$Dw\Big(D_{e_{a}}(D_{e_{b}}(r\nu
))\Big)=(1-\mu
_{b})Dw\Big(D_{e_{a}}e_{b}\Big),$$
$$\displaystyle Dw\Big(D_{D_{e_{a}}e_{b}}(r\nu )\Big)=(1-\mu
_{b})Dw\Big(D_{e_{a}}e_{b}\Big)$$
and
$$D_{ab}(rD_{\nu }w)=D^{3}w(e_{a},e_{b},r\nu )+(2-\mu _{a}-\mu _{b})D^{2}w(e_{a},e_{b})$$
which in view of the definition of the connexion $D$ and the definition of the covariant derivative implies that
\begin{equation}
D_{ab}(rD_{\nu}w)=rD_{\nu }(D_{ab}w)+(2-\mu _{a}-\mu _{b})D_{ab}w.
\end{equation}
In deriving this equality, we use the expression $(2.8)$ related to the curvature of $D$, in view of which the equality $\displaystyle D_{ab\nu }w=D_{\nu ab}w$ holds. 

\vskip2mm

Finally, we radially differentiate $(5.14)$ and multiply the result by $r$. Since
$D_{\nu }r=1$, $D_{\nu }u=D_{\nu}{\tilde w}=D_{\nu }v=0$ and taking $(5.16)$ and $(5.17)$ into account, we get
\begin{equation}
\begin{array}{cc}\displaystyle L[u](rD_{\nu }w_{s})+r\sum _{1\leq i,j\leq
n}G'^{ij}_{u}D_{ij}w_{s}-arD_{\nu}w_{s}=BrD_{\nu }w_{s}\\ \\ \displaystyle
+Br^2D_{\nu \nu}w_{s}+s\Big[2r\sum _{1\leq i,j\leq n}G'^{ij}_{u}D_{i}uD_{j}{\tilde
w_{s}}-arv\Big].\end{array}
\end{equation}
The equation of Gauss together with the computations of the second section show that for functions $w_{1},w_{2}\in \mathscr{C}^{2}({\tilde \Sigma })$ having the same values on $\Sigma $, we has
\begin{equation}
D_{ab}w_{1}=D_{ab}w_{2}+(1-\mu _{a})G_{ab}D_{\nu }(w_{1}-w_{2}),\ a,b\leq n+m-1,\mbox{ on }\Sigma .
\end{equation}
Since $w_{s}={\tilde w_{s}}$ and $D_{\nu }w_{s}=0$ everywhere on $\Sigma $, equation $(5.18)$ when restricted to $\Sigma$ implies the following :
\begin{equation}
\sum _{1\leq i,j\leq n}G'^{ij}_{u}(2sD_{i}uD_{j}{\tilde w_{s}}-D_{ij}{\tilde w_{s}})=sav.
\end{equation}
On the other hand, using $(5.19)$ and taking $(5.20)$ into account, we deduce from $(5.14)$ that everywhere on $\Sigma $, we have
\begin{equation}
\sum _{n+1\leq \alpha ,\beta \leq n+m-1}G'^{\alpha \beta }_{u}(2sD_{\alpha }uD_{\beta}{\tilde w_{s}}-D_{\alpha \beta }{\tilde w_{s}})+(a+\lambda){\tilde w_{s}}=sv.
\end{equation}
Since ${\tilde w_{s}}$ is a radially constant function, a combination of $(5.20)$ and $(5.21)$ shows that  ${\tilde w_{s}}$ is another solution of $(5.14)$. Hence, the maximum principle implies that ${\tilde w_{s}}=w_{s}$ everywhere on ${\tilde
\Sigma }$ and then $w_{s}$ is a radial constant.   

\vskip2mm

Now, restricting $(5.14)$ to $\Sigma $, using the maximum principle and the classical theory of uniformly elliptic equations [5], we establish the existence of a sequence of real numbers $R_{k}$ such that
\begin{equation}
\Vert w_{s}\Vert _{\mathscr{C}^{k}({\tilde \Sigma })}<R_{k},\ \mbox{for any }k\geq 0.
\end{equation}
Letting
$$\displaystyle B=\{ w\in \mathscr{C}^{\infty }(\tilde \Sigma )\mid \Vert
w\Vert _{\mathscr{C}^{k}(\tilde \Sigma  )}<R_{k},\mbox{ for any }k\geq 0\},$$
by virtue of $(5.22)$, equation $(5.13)$ has no solutions on the boundary of $B$ for any $s$. From Nagumo's theorem [7], the degree of $H$ at the origin with respect to
$B$ is invariant
\begin{equation}
d(H(s,.),0,B)=d(H(0,.),0,B)=\gamma .
\end{equation}
For $s=0$, $w_{0}=0$ is the unique solution of $(5.13)$ and standard computations show
that for $y\in \mathscr{C}^{\infty }(\tilde \Sigma )$, we have $d_{0}H_{0}y=w$, where $w$ is the unique solution of 
$$\left\{\begin{array}{cc}\displaystyle D_{\nu \nu}w+\sum _{a,b=1}^{n+m-1}G^{ab}D_{ab}w-(a+\lambda)w=D_{\nu \nu }y+r(m-1)D_{\nu }y\; \mbox{ in }{\tilde \Sigma }\\ \\ \displaystyle D_{\nu }w=0\; \mbox{ on }\partial {\tilde \Sigma }.\end{array}\right. $$ 
Arguing as above, we see that $\ker (\mbox{Id}-d_{0}H_{0})=\{0\}$. On the other hand, it is an easy fact to show that $\displaystyle d_{0}H_{0}$ is compact. Therefore, it follows from Fredholm's theorem and from the regularity of elliptic equations that $\mbox{Id}-d_{0}H_{0}$ is surjective. Thus, $\mbox{Id}-d_{0}H_{0}$ is invertible and $0$ is regular for $\mbox{Id}-H_{0}$.
Particularly, $\gamma =\pm 1$ in $(5.23)$ and $H_{1}$ has a fixed point $w$ solution in
$\Sigma $ of the following equations :
\begin{equation}
\left\{\begin{array}{cc}\displaystyle \sum _{1\leq i,j\leq n}G'^{ij}_{u}(2D_{i}uD_{j}w-D_{ij}w)=av\\ \\ \displaystyle \sum _{n+1\leq \alpha ,\beta \leq n+m-1}G'^{\alpha \beta
}_{u}(2D_{\alpha }uD_{\beta}w-D_{\alpha \beta }w)+(a+\lambda)w=v.\end{array}\right.
\end{equation}
Now, suppose $a\neq 0$ and choose a function $v$ which does not vanish identically on $\Sigma $. Then, for any $u\in \Theta $, the first equation in $(5.24)$ gives us
the existence of a function $w\in \mathscr{C}^{\infty }(\Sigma )$ such that
$$\sum _{1\leq i,j\leq n}G'^{ij}_{u}(2D_{i}uD_{j}w-D_{ij}w)\neq 0.$$
Hence $\Theta $ is a hypersurface.

\vskip2mm

To conclude, if $a=0$, we see from $(5.24)$ that $w\in B^{\infty }(\Sigma )$ solves $(5.10)$. Hence, for any $u\in \Theta $, $d_{u}\Gamma $ is invertible as a linear operator from $B^{\infty }(\Sigma )$ to $\mathscr{C}^{\infty }(\Sigma )$. We can then use in an obvious way the implicit function theorem of Nash-Moser [6] to show that $T$ is open in $[0,1]$. This ends the proof of the theorem. 

\vskip7mm

\subsection{Proof of theorem 4}

\vskip4mm

To prove this theorem, we apply the same topological argument as that used in the openness part of the proof of the previous theorem. Recall that we are given a strictly positive function $K\in \mathscr{C}^{\infty }(E_{*})$ and we suppose there exists two real numbers $r_1$ and $r_2$ such that $0<r_1\leq 1\leq r_2$ and  
\begin{equation}
\left\{ \begin{array}{c}\displaystyle K(\xi )>(\Vert \xi \Vert )^{-(m-1)} \mbox{ if
}\ \ \Vert \xi \Vert <r_{1}\\ \\ \displaystyle K(\xi )<(\Vert \xi \Vert )^{-(m-1)}\mbox{ if }\Vert \xi \Vert >r_{2}.\end{array}\right.
\end{equation}
Section 3 says that the problem we want to solve is equivalent to solving in $\mathscr{C}^{\infty}(\Sigma)$ the following equation
\begin{equation}
{\cal N}_{2}(u)=e^{(m-1)u}K(e^{u}\xi )(1+\vert D^{v}u\vert ^2)^{\frac{m+1}{2}}.
\end{equation}
For this purpose, to any $t\in [0,1]$ and any $w\in C^{\infty }(\Sigma )$, we associate
$H_{t}w=u_{t}$ where $u_{t}$ is the unique admissible solution of
\begin{equation}
\left\{ \begin{array}{c}\displaystyle {\cal N}_{1}(u)=1\\ \\ \displaystyle {\cal
N}_{2}(u)=e^{-u}\left[e^{mw}K(e^{w}\xi )\right]^{t}(1+\vert D^{v}u\vert ^2)^{\frac{m+1}{2}}.\end{array}\right.
\end{equation}
Let ${\cal B}$ be a bounded subset of $\mathscr{C}^{\infty }(\Sigma )$. Theorem 3 
ensures the existence of $u_{t}$ as well as that of a sequence of positive real numbers $C_{k}$ such that
\begin{equation}
\Vert u_{t}\Vert _{\mathscr{C}^{k}(\Sigma )}\leq C_{k}\ \mbox{for\ any\ }(t,w)\in [0,1]\times {\cal B}.
\end{equation}  
Thus, the operator $H$ defined by setting 
$$H(t,w)=H_{t}w\mbox{ for }(t,w)\in [0,1]\times \mathscr{C}^{\infty }(\Sigma )$$
is compact. Since $H_{0}(w)=0$ for all $w\in \mathscr{C}^{\infty }(\Sigma )$, the existence of a fixed point of $H(1,.)$, solution of equation $(5.26)$, reduces in
view of the theorem of Nagumo [7], to establishing that the set
$${\cal C}=\{u\in \mathscr{C}^{\infty }(\Sigma ): H(t,u)=u,\ t\in [0,1]\}$$
is bounded in $\mathscr{C}^{\infty }(\Sigma )$.

\vskip2mm

\noindent To deal with the $\mathscr{C}^0$-estimate, notice that any function $u\in {\cal C}$
satisfies
\begin{equation}
{\cal N}_{2}(u)=e^{(t-1)u}\left[e^{(m-1)u}K(e^{u}\xi )\right]^{t}(1+\vert D^{v}u\vert ^2)^{\frac{m+1}{2}}.
\end{equation}
So, at a point $\xi _{1}\in \Sigma $ where $u$ attains its maximum and since in a frame
adapted to $u$, we have $D_{\alpha }u(\xi _{1})=0$ and $D_{\alpha \alpha
}u(\xi _{1})\leq 0$, $(5.29)$ yields
$$1\leq e^{(t-1)u}\left[e^{(m-1)u}K(e^{u}\xi)\right]^{t}.$$
Thus, if $u\geq \log (r_2)$, by assumption $(5.25)$ we have
$$1<e^{(t-1)u}\leq (r_2)^{t-1}\leq 1,$$
which is a contradiction. Thus $u\leq u(\xi _1)\leq \log (r_2)$. The lower bound $u\geq \log (r_1)$ is obtained in a similar way.

\vskip2mm

\noindent The a priori estimates till order three are given in Lemmas 1, 2 and 3 and the higher order estimates may be established by induction as in the proof of the previous theorem. Consequently $H_{1}$ has a fixed point $u\in \mathscr{C}^{\infty }(\Sigma )$ which is an admissible solution of $(5.26)$. This completes the proof of the theorem.

\end{document}